\documentclass[twocolumn]{autart}    

\usepackage{graphicx}          
 \usepackage{lineno,hyperref}
\usepackage{amsmath}
\usepackage{amssymb}
\usepackage{amscd,multirow}
\usepackage{graphicx}
\usepackage{epstopdf}
\usepackage[titletoc]{appendix}
\usepackage{subfigure}
\usepackage[justification=centering]{caption}
\usepackage{cite,natbib}
\usepackage{wrapfig}
\usepackage{graphicx}
\usepackage{textcomp}
\usepackage{amsfonts}
\usepackage{subfigure}
\usepackage{subeqnarray,cases}
\usepackage{hyperref}
\usepackage[misc]{ifsym}
\usepackage{algorithm}
\usepackage{setspace}
\usepackage{xcolor}
\usepackage{subcaption}

\def\bigO{\mathcal{O}}

\hypersetup{hypertex=true,
            colorlinks=true,
            linkcolor=blue,
            anchorcolor=blue,
            citecolor=blue}
            
\allowdisplaybreaks[4]

\newtheorem{theorem}{Theorem}

\newtheorem{lemma}{Lemma}

\newtheorem{remark}{Remark}

\newtheorem{assumption}{Assumption}

\newenvironment{proof}{\noindent {\textbf{Proof.}}}

\begin{document}

\begin{frontmatter}

\title{Nesterov acceleration for strongly convex-strongly concave bilinear saddle point problems: discrete and continuous-time approaches\thanksref{footnoteinfo}} 

\thanks[footnoteinfo]{This work was supported by the Talent Introduction Project of Xihua University (Grant No. RX2400002043), Sichuan Science and Technology  Program (Grant No. 2025ZNSFSC0813) and the National Natural Science Foundation of China (Grant No. 12471296). Corresponding author: Ya-Ping Fang.}

\author[1]{Xin He}\ead{hexinuser@163.com},
\author[2]{Ya-Ping Fang}\ead{ypfang@scu.edu.cn} 

\address[1]{School of Science, Xihua University, Chengdu, Sichuan, 610039,  China}
\address[2]{Department of Mathematics, Sichuan University, Chengdu, Sichuan, 610064, China}



\begin{abstract}    
{In this paper, we study a bilinear saddle point problem of the form
$\min_{x}\max_{y} F(x) + \langle Ax, y \rangle - G(y)$, where $F$ and $G$ are $\mu_F$- and $\mu_G$-strongly convex functions, respectively. By incorporating Nesterov acceleration for strongly convex optimization, we first propose a  first-order  primal-dual gradient algorithm. We show that it achieves the linear convergence rate $\mathcal{O}((1 - \min\{\sqrt{{\mu_F}/{L_F}}, \sqrt{{\mu_G}/{L_G}}\})^k)$ for both the primal-dual gap and the iterative gap, where $L_F$ and $L_G$ denote the smoothness constants of $F$ and $G$, respectively. We further develop a continuous-time accelerated primal-dual dynamical system with constant damping. Using the Lyapunov analysis method, we establish the linear convergence rate $\mathcal{O}(e^{-\min\{\sqrt{\mu_F},\sqrt{\mu_G}\}t})$. Notably, when $A = 0$, our methods recover the classical Nesterov accelerated methods for strongly convex unconstrained problems in both discrete and continuous-time. Numerical experiments are presented to support the theoretical convergence results.}
\end{abstract}

\end{frontmatter}

\section{Introduction}

In this paper, we consider the strongly convex-strongly concave saddle point problem with bilinear coupling, given by: 
\begin{equation}\label{ques} 
	\min_{x\in\mathbb{R}^n}\max_{y\in\mathbb{R}^m} \mathcal{L}(x,y) = F(x) + \langle Ax, y \rangle - G(y), 
\end{equation} 
where $F$ and $G$ are strongly convex functions, and $A \in \mathbb{R}^{m \times n}$ is a coupling matrix. Problem \eqref{ques} appears in a number of applications involving bilinear coupling between primal and dual variables, including policy evaluation in reinforcement learning, regularized learning models, and quadratic minimax problems in control and numerical optimization \citep{CherukuriSICON,Du2017,Kovalev2022,ZhaoJMLR,ZhangICAIS}. This motivates the study of efficient primal-dual methods for solving problem \eqref{ques}. Our discussion begins with three representative examples.

{\bf Reinforcement learning.} 
In reinforcement learning, policy evaluation often involves minimizing the mean squared projected Bellman error (MSPBE)~\citep{Du2017}. 
Given samples $\{s_t,a_t,r_t,s_{t+1}\}_{t=1}^n$ under a policy $\pi$, define 
$\mathbf{A}=\tfrac{1}{n}\sum_t\boldsymbol{\phi}(s_t)(\boldsymbol{\phi}(s_t)-\gamma\boldsymbol{\phi}(s_{t+1}))^T$, 
$\mathbf{b}=\tfrac{1}{n}\sum_t r_t\boldsymbol{\phi}(s_t)$, and 
$\mathbf{C}=\tfrac{1}{n}\sum_t\boldsymbol{\phi}(s_t)\boldsymbol{\phi}(s_t)^T$, 
where $\boldsymbol{\phi}(s_t)$ is the feature vector and $\gamma$ the discount factor.  Then, the empirical MSPBE minimization can be expressed as  
\[
\min_{\boldsymbol{\theta}}\max_{\mathbf{w}}
\Big(\tfrac{\gamma}{2}\|\boldsymbol{\theta}\|^2
-\mathbf{w}^T\mathbf{A}\boldsymbol{\theta}
-\tfrac{1}{2}\|\mathbf{w}\|_{\mathbf{C}}^2
+\mathbf{w}^T\mathbf{b}\Big),
\]
which is a strongly convex-strongly concave instance of~\eqref{ques} when $\mathbf{C}$ is positive definite.

{\bf Quadratic minimax problem.}
 Quadratic minimax problems arise in several contexts, including numerical analysis and optimal control \citep{Wang2020, LiuTNNLS}. Consider the quadratic forms $F(x) = x^T R x$ and $G(y) = y^T S y$ with $R\succ 0$ and $S\succ0$ being positive definite matrices. The resulting minimax objective is quadratic in both $x$ and $y$ and is given by
\begin{equation}\label{ques_quad}
\min_x\max_y\mathcal{L}(x,y) = x^TRx + \langle  Ax, y\rangle  - y^TSy.
\end{equation}
Although \eqref{ques_quad} is quadratic in both variables, it is still nontrivial to solve efficiently. The main reason is that the bilinear coupling term $\langle Ax,y\rangle$ creates a nonseparable interaction between the primal and dual variables, so the problem cannot be reduced to two independent strongly convex minimization problems. This coupling  makes the design and analysis of accelerated primal-dual methods substantially more delicate \citep{ZhangICAIS}.

{\bf Empirical risk minimization.}
Problem \eqref{ques} has an important application in empirical risk minimization (ERM), a central model in machine learning \citep{Shalev2014}. A standard ERM problem takes the form
\[    \min_{x \in \mathbb{R}^{n}} F(x) + H(Ax),\]
where $H$ is a convex loss function, $A$ is the data matrix, and $F$ is a strongly convex regularizer. By Fenchel duality, this problem can be equivalently reformulated as  saddle-point problem \eqref{ques}:
\[
    \min_{x}\max_{y} F(x) + \langle Ax,y\rangle - H^*(y),
\]
where $H^*$ is the Fenchel conjugate of $H$.

For the unconstrained strongly convex optimization problem $\min_x F(x)$, where $F$ is $\mu$-strongly convex and $L$-smooth (i.e., $\nabla F$ is $L$-Lipschitz continuous),  \cite{Nesterov2003} proposed the Nesterov accelerated gradient method (NAG-\texttt{SC}):
\begin{equation}\label{al_Nag_sc}
\begin{aligned}
     & \bar{x}_{k} = x_{k} + \tfrac{1 - \sqrt{\mu r} }{ 1 + \sqrt{\mu  r} } \left( x_{k} - x_{k-1} \right),\\
         & x_{k + 1} =  \bar{x}_{k} - r \nabla F(\bar{x}_{k}),
         \end{aligned}
\end{equation}
where $r\leq 1/L$. It shows a linear convergence rate $F(x_k)-\min F= O((1-\sqrt{\mu r})^k)$. When $r=1/L$, this method achieves the  linear convergence rate $O((1-\sqrt{\mu/L})^k)$, improving on the classical gradient descent method \citep{Nesterov2003,LinML2019}. Taking the limit $r \rightarrow 0$, one obtains the continuous-time dynamic counterpart \citep{Wilson}: 
\begin{equation}\label{dy_low_hb}
\ddot{x}(t) + 2 \sqrt{\mu}\dot{x}(t) +  \nabla F(x(t)) = 0.
\end{equation}      
\citet{LuoMP}, and  \citet{Wilson}, established that this dynamic exhibits linear convergence $F(x(t)) - \min F \leq O(e^{-\sqrt{\mu}t}).$

For saddle point problem \eqref{ques},  $F$ and $G$ satisfy the following assumption:
\begin{assumption}\label{ass1}
$F:\mathbb{R}^n\to \mathbb{R}$ is $\mu_F$-strongly convex and $L_F$-smooth, $G:\mathbb{R}^m\to \mathbb{R}$ is $\mu_G$-strongly convex and $L_G$-smooth.
\end{assumption}

Under Assumption \ref{ass1}, although various gradient-based methods for solving \eqref{ques} are known to enjoy linear convergence rates depending on the strong convexity parameters, no Nesterov-type accelerated primal-dual gradient method has yet been proved to attain the rate
\begin{equation}\label{opt_conv}
	\bigO\left(\left(1-\min\left\{\sqrt{\tfrac{\mu_F}{L_F}},\sqrt{\tfrac{\mu_G}{L_G}}\right\}\right)^k\right)
\end{equation}
for the primal-dual gap $\mathcal{L}(x_k,y^*)-\mathcal{L}(x^*,y_k)$. This rate further yields the oracle complexity
\begin{equation}\label{opt_oracle}
\bigO\left(\max\left\{\sqrt{\tfrac{L_F}{\mu_F}},\sqrt{\tfrac{L_G}{\mu_G}}\right\}\log\tfrac{1}{\varepsilon}\right)
\end{equation}
in terms of gradient evaluations of $\nabla F$ and $\nabla G$ for computing an $\varepsilon$-solution. Table \ref{tab1} compares this complexity bound with the oracle complexity bounds of several existing methods for solving problem~\eqref{ques} under Assumption~\ref{ass1}.
\begin{table}[t]
\caption{Comparison of  oracle complexity bounds}
\label{tab1}
\centering
\renewcommand{\arraystretch}{1.2}
\begin{tabular}{p{0.95\columnwidth}}
\hline
\textbf{Reference and complexity} \\
\hline
\cite{Cohen2021}:  \\
\quad $\bigO\!\left(\left(\frac{L_F}{\mu_F}+\frac{\|A\|}{\sqrt{\mu_F\mu_G}}+\frac{L_G}{\mu_G}\right)\log\frac{1}{\varepsilon}\right)$ \\
\hline

\cite{Wang2020}:  \\
\quad $\widetilde{\bigO}\!\left(\left(\sqrt{\frac{L_F}{\mu_F}}+\sqrt{\frac{\|A\|\max\{L_F,\|A\|,L_G\}}{\mu_F\mu_G}}+\sqrt{\frac{L_G}{\mu_G}}\right)\log\frac{1}{\varepsilon}\right)$ \\
\hline

\cite{Xie2020}:  \\
\quad $\widetilde{\bigO}\!\left(\left(\frac{(L_F^2L_G)^{1/4}}{\mu_F^{1/2}\mu_G^{1/4}}+\frac{\|A\|}{\sqrt{\mu_F\mu_G}}+\frac{(L_G^2L_F)^{1/4}}{\mu_G^{1/2}\mu_F^{1/4}}\right)\log\frac{1}{\varepsilon}\right)$ \\
\hline

\cite{Thekumparampil2022}, \cite{Kovalev2022}, \cite{Du2022}:  \\
\quad $\bigO\!\left(\left(\sqrt{\frac{L_F}{\mu_F}}+\frac{\|A\|}{\sqrt{\mu_F\mu_G}}+\sqrt{\frac{L_G}{\mu_G}}\right)\log\frac{1}{\varepsilon}\right)$ \\
\hline

This paper: $\bigO\!\left(\max\!\left\{\sqrt{\frac{L_F}{\mu_F}},\sqrt{\frac{L_G}{\mu_G}}\right\}\log\frac{1}{\varepsilon}\right)$ \\
\hline
\end{tabular}
\end{table}

We first demonstrate that under Assumption~\ref{ass1}, the linear convergence rate ~\eqref{opt_conv} of a first-order primal-dual gradient algorithm  to solve problem~\eqref{ques} may achieved. Consider the case $A=0$ in  problem \eqref{ques}, it reduces to two independent strongly convex optimization problems $\min_x F(x)$ and $\min_y G(y)$. Applying  NAG-\texttt{SC} \eqref{al_Nag_sc} (with $r=1/L$) to solve separately $\min_x F(x)$ and $\min_y G(y)$ yields the linear convergence rate~\eqref{opt_conv} for the primal-dual gap $\mathcal{L}(x_k,y^*)-\mathcal{L}(x^*,y_k)$. This naturally raises the question: Can the acceleration technique of NAG-\texttt{SC} be extended to strongly convex-strongly concave saddle point problems \eqref{ques} with the linear convergence rate \eqref{opt_conv}? In this work, we explore this question in depth.

 In recent years, continuous-time inertial dynamical systems have been studied on the design and analysis of various damping~\citep{BoSICON,AttouchJEMS,SuJMLR,HeNN,HeTAC}.
Compared with a discrete algorithm, a continuous-time dynamical system serves primarily as a theoretical framework that uncovers the fundamental convergence mechanism underlying algorithmic acceleration. It clarifies the roles of inertial and damping terms in determining stability and convergence rates, and provides a continuous counterpart that helps interpret and justify the acceleration effects observed in the discrete-time algorithm.  Moreover, insights gained from a continuous-time model often guide the design of new and faster discrete algorithms. Building on Nesterov's foundational work on accelerated methods for strongly convex optimization~\citep{Nesterov2003}, a line of research has focused on inertial dynamical systems with constant viscous damping to solve unconstrained strongly convex problems~\citep{ShiNIPS,LuoMP,Wibisono2016,Siegel}. In particular, many inertial primal-dual dynamical systems have been developed to solve the linearly constrained problem:
\begin{equation}\label{ques_cons}
    \min_{x} F(x), \quad \text{s.t.} \quad Ax = b,
\end{equation}
which can be interpreted as a special case of the general saddle point problem~\eqref{ques} by setting $G(y) = \langle b, y \rangle$. Most existing second-order systems address convex objectives and achieve, at best, a convergence rate of $\mathcal{O}(1/t^2)$ without time-scaling coefficients; see~\citet{ZengTAC,HeSICON,BotJDE2022,HeNN,HeAuto}. \citet{HeAuto2025} were the first to introduce a continuous-time ``second-order primal + first-order dual'' dynamic for solving~\eqref{ques_cons} in the strongly convex case, proving $\mathcal{O}(e^{-\sqrt{\mu}t})$ linear convergence. Meanwhile, considerable efforts have been made to extend these inertial dynamics from the constrained setting~\eqref{ques_cons} to the general bilinear saddle point problem~\eqref{ques} in convex case \citep{ZengIFAC,DingCoa,HeAMO,Luo2024Arx,LiZeng2024,SunXK}.  In contrast, the accelerated scheme NAG-\texttt{SC} plays a fundamental role in the strongly convex setting, laying the foundation for both discrete and continuous-time methods. Despite its importance, accelerated primal-dual continuous-time dynamics for strongly convex-strongly concave saddle point problems remain largely unexplored. In this paper, we propose a new continuous-time inertial primal-dual dynamical system for solving~\eqref{ques} in the strongly convex-strongly concave setting.

 Although problem~\eqref{ques} is related to broader equilibrium models, including Nash equilibrium seeking, our objective here is more focused. We study the strongly convex-strongly concave bilinear saddle-point problem as a basic primal-dual model, and investigate whether Nesterov-type acceleration can be established with explicit accelerated linear convergence guarantees in both discrete time and continuous time. The present work is developed in a centralized and full-information setting.

{\bf Main contributions.} Compared with unconstrained strongly convex optimization or linearly constrained minimization, the main difficulty here lies in obtaining Nesterov-type accelerated rates under the coupled primal-dual structure induced by the bilinear term. The interaction between $x$ and $y$ through the operator $A$ makes it difficult to construct a single energy function that simultaneously controls both variables and yields a decay rate comparable to NAG-\texttt{SC}. The proposed discrete algorithm and continuous-time dynamic are designed to address this issue by balancing the primal-dual coupling in both the iteration scheme and the Lyapunov analysis. Our main contributions are summarized as follows:

\begin{itemize}
\item[(a)] \textbf{Discrete acceleration:}
Building on the Nesterov acceleration framework, we propose a first-order primal-dual gradient method for solving the bilinear saddle point problem~\eqref{ques}, where the functions $F$ and $G$ satisfy Assumption~\ref{ass1}. A key feature of the method is that it combines a symmetric momentum construction for both variables with a shared extrapolation parameter and an additional coupling extrapolation term. This design makes it possible to balance the interaction between the primal and dual updates and to establish a single contractive energy estimate for the coupled system. We show that the proposed method achieves the linear convergence rate~\eqref{opt_conv} for both the primal-dual gap and the iterate error, thereby improving upon the linear convergence guarantees of existing first-order methods~\citep{Borodich2023,Cohen2021,Kovalev2022,Thekumparampil2022,Du2022}. 

\item[(b)] \textbf{Continuous-time acceleration:}
We further propose a continuous-time inertial primal-dual dynamical system tailored to problem~\eqref{ques}. The dynamic is designed so that the damping and coupling terms can be incorporated into a unified Lyapunov analysis, allowing the primal-dual interaction to be controlled in a way analogous to accelerated dynamics for strongly convex minimization. We establish linear convergence at the rate $\mathcal{O}\big(e^{-\min\{\sqrt{\mu_F}, \sqrt{\mu_G}\}t}\big)$ for the primal-dual gap, the trajectory error, and the velocity norm. This result extends the inertial dynamics in~\citet{Siegel,LuoMP,HeAuto2025} to the saddle-point setting. Moreover, the proposed continuous-time dynamic provides a theoretical framework for interpreting the acceleration mechanism of the discrete algorithm. 
\end{itemize}

{\bf Notations:} Let \( \langle \cdot, \cdot \rangle \) and \( \|\cdot\| \) denote the inner product and the Euclidean norm, respectively. \( \mathbf{I}_n \) represents the \( n \times n \) identity matrix. For a symmetric positive definite matrix $R$, we denote  $\lambda_{\max}(R)$ and $\lambda_{\min}(R)$ be the largest and smallest eigenvalues of $R$, respectively; $\kappa(R)= \tfrac{\lambda_{\max}(R)}{\lambda_{\min}(R)}$ be a condition number of $R$. For a differentiable function \( F: \mathbb{R}^{n} \to \mathbb{R} \), we say that \( F \) is \( \mu_F \)-strongly convex if for any  $x, y \in \mathbb{R}^n$, $F(y) - F(x) - \langle \nabla F(x), y - x \rangle \geq \frac{\mu_F}{2} \|y - x\|^2$, and \( F \) is \( L_F \)-smooth if $
\| \nabla F(y) - \nabla  F(x) \| \leq L_F \| y - x \|$.
Since \( F \) and \( G \) are strongly convex, then problem \eqref{ques} admits a unique solution, which we denote as \( (x^*, y^*) \in \mathbb{R}^n \times \mathbb{R}^m \), such that for any $(x,y)\in \mathbb{R}^n\times  \mathbb{R}^m$,
\begin{equation}\label{eq_lag}
	\mathcal{L}(x^*,y)\leq \mathcal{L}(x^*,y^*)\leq \mathcal{L}(x,y^*).	
\end{equation}

 {\bf Organization:} Section~\ref{sec2} introduces a  first-order primal-dual gradient method for solving  problem~\eqref{ques} and establishes its convergence rate. Section~\ref{sec3} presents the  continuous-time dynamical system along with the  theoretical results. Section~\ref{sec4} provides numerical experiments to demonstrate the effectiveness of the proposed methods. Finally, Section~\ref{sec5} concludes the paper.

\section{Nesterov accelerated  first-order primal-dual gradient algorithm}\label{sec2}

In this section,  we always assume that Assumption \ref{ass1} holds.  By examining  NAG-\texttt{SC} \eqref{al_Nag_sc} for unconstrained strongly convex optimization, we propose a first-order primal-dual gradient method for solving the smooth strongly convex-strongly concave problem \eqref{ques}. Specifically, we consider the iterative sequence  $\{(x_k, y_k)\}_{k \geq 1}$ that satisfies the following equations:
\begin{small}
\vspace{-1em}
\begin{subequations}\label{eq_method}
\begin{numcases}{}
	(\bar{x}_k, \bar{y}_k)= (x_k, y_k)  + \tfrac{1 - \theta}{1 + \theta} \left[ (x_k, y_k) - (x_{k-1}, y_{k-1}) \right], \label{eq_method_a} \\
	x_{k+1} = \bar{x}_k  - r \left( \nabla F(\bar{x}_k) + A^T \left( y_k + \tfrac{1}{\theta}(y_{k+1} - y_k) \right) \right), \label{eq_method_b} \\
	y_{k+1} = \bar{y}_k - s\left( \nabla G(\bar{y}_k) - A \left( x_k + \tfrac{1}{\theta}(x_{k+1} - x_k) \right) \right), \label{eq_method_c}
\end{numcases}
\end{subequations}
\end{small}
where
\[ \theta = \min\left\{ \sqrt{\mu_F r}, \sqrt{\mu_G s } \right\}\]
with $r\leq 1/L_F$ and $s\leq 1/L_G$. Under Assumption \ref{ass1}, it is easy to verify that $\theta\in(0,1]$, see \cite{Nesterov2003}. Comparing equation \eqref{eq_method} with NAG-\texttt{SC}, we observe the following two key differences:

1) \textbf{Interpolation in the momentum term:} In NAG-\texttt{SC}, when \( \bar{x}_k \) is interpolated using \( x_k \) and \( x_{k-1} \), the interpolation coefficient is \( \tfrac{1 - \theta}{1 + \theta} \), with \( \theta = \sqrt{\mu_F r} \) and $r\leq 1/L$. However, in equation \eqref{eq_method}, we define $\theta = \min\left\{ \sqrt{\mu_F r}, \sqrt{\mu_G s } \right\}$ with $r\leq 1/L_F$ and $s\leq 1/L_G$. This choice is due to the symmetry of the functions \( F \) and \( G \) in the saddle point problem. Therefore, when selecting the interpolation coefficients, the characteristics of both functions must be considered simultaneously. When \( A = 0 \),  the \( F \)-part of equation \eqref{eq_method_a} corresponds to NAG-\texttt{SC}. Similar interpolation techniques for strongly convex optimization are discussed in \cite{LinML2019, HeAuto2025, ShiMP}.

2)  \textbf{Extra extrapolation term:} Unlike the first-order primal-dual gradient methods in \cite{Kovalev2022, Thekumparampil2022, Borodich2023}, we include the extrapolation term \( \tfrac{1}{\theta}(x_{k+1} - x_k) \) and \( \tfrac{1}{\theta}(y_{k+1} - y_k) \) in equations \eqref{eq_method_b} and \eqref{eq_method_c}, respectively. This term plays a crucial role in the convergence analysis. Note that the extrapolation term has been widely used in designing accelerated primal-dual gradient methods for solving linearly constrained problems in the convex case (see \cite{HeAuto, LuoJota, BotMP,LuoMC,HeANM}).

By examining equation \eqref{eq_method}, we see that the sequences $x_{k+1}$ and $y_{k+1}$ are coupled with each other and cannot yield a useful iterative sequence. However, we can substitute $y_{k+1}$ from equation \eqref{eq_method_c} into equation \eqref{eq_method_b} to obtain 
\begin{eqnarray*}
	&&(\mathbf{I}_n + \tfrac{rs}{\theta^2} A^T A) x_{k+1}\\
	&&\quad= \bar{x}_k - r \left( \nabla F(\bar{x}_k) + A^T \left( \left( 1 - \tfrac{1}{\theta} \right) y_k + \tfrac{1}{\theta} \hat{y}_k \right) \right),
\end{eqnarray*}
 where $\hat{y}_k = \bar{y}_k - s \left( \nabla G(\bar{y}_k) - \left( 1 - \tfrac{1}{\theta} \right) A x_k \right)$.
Based on this, we propose a accelerated primal-dual gradient method (Algorithm \ref{al_grad}). It is easy to verify that the iterative sequence of Algorithm \ref{al_grad} also satisfies equation \eqref{eq_method}.
\begin{algorithm}[h] 
\caption{Nesterov Accelerated First-Order Primal-Dual Gradient Method}
\label{al_grad}
{\bf Input:} $x_1=x_0 \in \mathbb{R}^{n}$, $y_1=y_0\in \mathbb{R}^{m}$, $\mu_F>0$, $L_F>0$, $\mu_G>0$, $L_G>0$, $A\in \mathbb{R}^{m\times n}$.

\vspace{0.1em}
{\bf Parameters:} Set $r\leq 1/L_F$, $s\leq 1/L_G$ and
\[
\theta = \min\left\{ \sqrt{\mu_F r}, \sqrt{\mu_G s } \right\},\quad 
\mathcal{B} = \left(\mathbf{I}_n+\tfrac{rs}{\theta^2}A^TA\right)^{-1}.
\]

{\bf For} {$k = 1, 2,\cdots, K-1$} {\bf do}
\[
\left\{
\begin{array}{l}
(\bar{x}_k,\bar{y}_k) 
= (x_k,y_k)+\tfrac{1-\theta}{1+\theta}\left[(x_k,y_k)-(x_{k-1},y_{k-1})\right]\\[0.2em]
\hat{y}_k
= \bar{y}_k-s\left(\nabla G(\bar{y}_k)-\left(1-\tfrac{1}{\theta}\right)Ax_k\right)\\[0.2em]
x_{k+1}
= \mathcal{B}\left(\bar{x}_k-r\left(\nabla F(\bar{x}_k)+A^T\left(\left(1-\tfrac{1}{\theta}\right)y_k+\tfrac{1}{\theta}\hat{y}_k\right)\right)\right)\\[0.2em]
y_{k+1}
= \bar{y}_k-s\left(\nabla G(\bar{y}_k)-A\left(x_k+\tfrac{1}{\theta}(x_{k+1}-x_k)\right)\right).
\end{array}
\right.
\]
\vspace{-1em}

{\bf end}

\vspace{0.1em}
{\bf Output:} $x_K$, $y_K$.
\end{algorithm}

To establish the fast convergence property of Algorithm~\ref{al_grad}, we outline a proof based on the energy analysis technique. 
The key idea is to construct a positive and nonincreasing energy sequence $\{\mathcal{E}_k\}_{k\geq 1}$ defined as
\begin{equation}\label{eq_ener_dis}
	\mathcal{E}_k = \mathcal{L}(x_k,y^*) - \mathcal{L}(x^*,y_k) + \tfrac{1}{2r}\|u_k\|^2 + \tfrac{1}{2s}\|v_k\|^2,
\end{equation}
where  the sequence \( \{(x_k, y_k)\}_{k \geq 1} \) is generated by Algorithm~\ref{al_grad}, \( (x^*, y^*) \) is a unique solution of problem \eqref{ques},  and
\begin{equation}\label{eq_uvk}
	\begin{aligned}
		u_k &= \theta(x_k - x^*) + (1 - \theta)(x_k - x_{k-1}), \\
		v_k &= \theta(y_k - y^*) + (1 - \theta)(y_k - y_{k-1}).
	\end{aligned}
\end{equation}
The construction of $\{\mathcal{E}_k\}_{k\geq 1}$ is inspired by the potential sequences employed in accelerated primal-dual methods for the linearly constrained convex problem~\eqref{ques_cons}~\citep{BotMP,HeNA}, which can be viewed as discrete counterparts of the continuous-time energy functions studied in~\citep{HeSICON,BotJDE2022}. 
However, the essential difference lies in the structural complexity of  saddle-point problem~\eqref{ques}. 
Unlike the linearly constrained convex case, where the primal variable evolves independently within a single convex objective, our setting involves two interdependent variables that interact through the bilinear term $\langle Ax, y\rangle$ under the strongly convex-strongly concave regime. 
To capture this coupling effect, the proposed energy sequence~\eqref{eq_ener_dis} incorporates a constant extrapolation parameter $\theta$ and introduces  joint scaling coefficients $\tfrac{1}{2r}$ and $\tfrac{1}{2s}$, to balance the primal and dual residuals. 
By exploiting the strong convexity and smoothness of $F$ and $G$, together with the iterative relations in~\eqref{eq_method} and an appropriate scaling argument, we prove that	$\mathcal{E}_{k+1} \le (1 - \theta)\mathcal{E}_k$.
This contraction inequality directly leads to the linear convergence rate of Algorithm~\ref{al_grad}.

We begin by demonstrating an inequality for smooth strongly convex smooth functions.

\begin{lemma}\label{lemmaSF}
Assume that \( F: \mathbb{R}^{n} \to \mathbb{R} \) is \( \mu_F \)-strongly convex and \( L_F\)-smooth. Then, the following inequality is valid for any $x,y,z\in\mathbb{R}^n$:
			\[F(y)-F(x)\leq\langle\nabla F(z), y-x\rangle+\tfrac{L_F}{2}\|y-z\|^2-\tfrac{\mu_F}{2}\|z-x\|^2.\]			
\end{lemma}
\begin{proof}
Since \( F \) is  \( \mu_F \)-strongly convex and \( L_F\)-smooth, as shown in \cite{Nesterov2003}, we can get that
\[
\tfrac{\mu_F}{2}\|x-y\|^2\leq F(y) - F(x)-\langle \nabla F(x),y-x\rangle \leq \tfrac{L_F}{2}\|x-y\|^2
\]
holds for any $x,y\in\mathbb{R}^n$.
This leads to 
\begin{eqnarray*}
   &&  F(y)-F(x) = F(y)-F(z)+F(z)-F(x)\\
	&&\quad \leq \langle \nabla F(z), y-z\rangle+\tfrac{L_F}{2}\|y-z\|^2\\
	&&\qquad+ \langle \nabla F(z), z-x\rangle-\tfrac{\mu_F}{2}\|z-x\|^2\\
	&&\quad 	=\langle\nabla F(z), y-x\rangle+\tfrac{L_F}{2}\|y-z\|^2-\tfrac{\mu_F}{2}\|z-x\|^2,
\end{eqnarray*}
and then the result holds.\hspace*{\fill} \textbf{$\square$}
\end{proof}

We now turn our attention to analyzing the convergence properties of Algorithm \ref{al_grad}.

\begin{lemma}\label{lemmaL}
Let \( \{(x_k, y_k)\}_{k \geq 1} \) be the sequence generated by Algorithm \ref{al_grad}. Define the energy sequence $\{\mathcal{E}_k\}_{k\geq 1}$ as in \eqref{eq_ener_dis}.  Then, for any $k\geq 1$, we have 
\[\mathcal{E}_{k+1} \leq  (1- \theta)\mathcal{E}_{k}.\]	
\end{lemma}
\begin{proof}
Since \( F \) is \( \mu_F \)-strongly convex and \( L_F \)-smooth, by Lemma \ref{lemmaSF} we get that for any $ x \in \mathbb{R}^n$,
  \begin{eqnarray}\label{eqH1}
 	F(x_{k+1})-F(x)&\leq&\langle\nabla F(\bar{x}_k), x_{k+1}-x\rangle \\
 	&&+\tfrac{L_F}{2}\|x_{k+1}-\bar{x}_k\|^2-\tfrac{\mu_F}{2}\|\bar{x}_k-x\|^2.\nonumber
 \end{eqnarray}
By the expression in equation \eqref{eq_method_b}, we can get
\[\nabla F(\bar{x}_k) = -\tfrac{1}{r}(x_{k+1}-\bar{x}_k)-A^T\left(y_k+\tfrac{1}{\theta}(y_{k+1}-y_k)\right).\]
This together with \eqref{eqH1} implies
\begin{eqnarray*}
&&	 F(x_{k+1})-F(x)\leq -\tfrac{1}{r}\langle x_{k+1}-\bar{x}_k, x_{k+1}-x\rangle\nonumber \\
	&&\qquad-\left\langle x_{k+1}-x, A^T\left(y_k+\tfrac{1}{\theta}(y_{k+1}-y_k)\right)\right\rangle\\
	&&\qquad+\tfrac{L_F}{2}\|x_{k+1}-\bar{x}_k\|^2-\tfrac{\mu_F}{2}\|\bar{x}_k-x\|^2\\
	&&\quad=-\tfrac{1}{r}\langle x_{k+1}-\bar{x}_k, \bar{x}_k-x\rangle\\
	&&\qquad-\left\langle x_{k+1}-x, A^T\left(y_k+\tfrac{1}{\theta}(y_{k+1}-y_k)\right)\right\rangle\nonumber \\
	&&\qquad-\left(\tfrac{1}{r}-\tfrac{L_F}{2}\right)\|x_{k+1}-\bar{x}_k\|^2-\tfrac{\mu_F}{2}\|\bar{x}_k-x\|^2.\nonumber
\end{eqnarray*}
By adding \( \theta \) times the above inequality with \( x = x^* \) and \( (1 - \theta) \) times the same inequality with \( x = x_k \), we obtain
\begin{eqnarray}\label{eqH2}
	&& F(x_{k+1})-F(x^*)-(1-\theta)(F(x_{k})-F(x^*)) \nonumber \\
	 && = \theta (F(x_{k+1})-F(x^*)) +(1-\theta)(F(x_{k+1})-F(x_{k}))\nonumber\\
	  && \leq -\tfrac{1}{r}\langle x_{k+1}-\bar{x}_k, \bar{x}_k-\theta x^*-(1-\theta) x_k\rangle \nonumber\\
	   &&  -\left\langle x_{k+1}-x_k+\theta(x_k-x^*), A^T(y_k+\tfrac{1}{\theta}(y_{k+1}-y_k))\right\rangle \nonumber\\
	   &&  - \tfrac{1}{2r}\|x_{k+1}-\bar{x}_k\|^2-\tfrac{\mu_F\theta}{2}\|\bar{x}_k-x^*\|^2\\
	   &&-\tfrac{\mu_F(1-\theta)}{2}\|\bar{x}_k-x_k\|^2,\nonumber
\end{eqnarray}
where the last equality follows from $r\leq 1/L_F$. 

Next, starting from \eqref{eq_method_a} and \eqref{eq_uvk}, we derive
\begin{eqnarray*}
	&& \tfrac{1}{2r}\| \bar{x}_{k}-\theta x^*-(1-\theta)x_k\|^2 \\
	&&= \tfrac{\theta^2}{2r}\left\| \bar{x}_{k}-x^*+\tfrac{1-\theta}{\theta}(\bar{x}_{k}-x_k)\right\|^2\\
	 && = \tfrac{\theta^2}{2r}\left\|\theta(\bar{x}_{k}-x^*)+(1-\theta)\left(\tfrac{1+\theta}{\theta}\bar{x}_{k}  - x^*-\tfrac{1}{\theta}x_k\right)\right\|^2\\
	  && = \tfrac{\theta^2}{2r}\left\|\theta(\bar{x}_{k}-x^*)+(1-\theta)\tfrac{u_k}{\theta}\right\|^2\\
	  &&\leq \tfrac{\theta^3}{2r}\|\bar{x}_{k}-x^*\|^2+\tfrac{(1- \theta)}{2r}\|u_k\|^2,
\end{eqnarray*}
where the last inequality follows from the convexity of 
 $\|\cdot\|^2$ and the fact that $\theta\in(0,1]$. Combining this with \eqref{eq_uvk}, we obtain
\begin{eqnarray*}
&&	-\tfrac{1}{r}\langle x_{k+1}-\bar{x}_k, \bar{x}_k-\theta x^*-(1-\theta) x_k\rangle \nonumber \\
	&&\qquad\quad = \tfrac{1}{2r} \|x_{k+1}-\bar{x}_k\|^2-\tfrac{1}{2r}\|u_{k+1}\|^2\nonumber \\
	&&\qquad\qquad+\tfrac{1}{2r}\|\bar{x}_k-\theta x^*-(1-\theta) x_k\|^2\nonumber \\
	&&\qquad\quad \leq  \tfrac{1}{2r}\|x_{k+1}-\bar{x}_k\|^2+\tfrac{\theta^3}{2r}\|\bar{x}_{k}-x^*\|^2\\
	&&\qquad\qquad+\tfrac{(1- \theta)}{2r}\|u_k\|^2- \tfrac{1}{2r}\|u_{k+1}\|^2,\nonumber
\end{eqnarray*}
where the first equality in \eqref{eq_24_1} is derived from
\[-\langle x,y\rangle = \tfrac{1}{2}(\|x\|^2+\|y\|^2-\|x+y\|^2),\quad \forall x,y\in \mathbb{R}^n.\] 
Then, from \eqref{eqH2} and $\theta\leq 1$, we have
\begin{eqnarray}\label{eq_24_1}
	&& F(x_{k+1})-F(x^*)-(1-\theta)(F(x_{k})-F(x^*)) \nonumber \\
		  && \leq \tfrac{\theta(\theta^2-\mu_Fr)}{2r}\|\bar{x}_k-x^*\|^2+\tfrac{(1- \theta)}{2r}\|u_k\|^2- \tfrac{1}{2r}\|u_{k+1}\|^2 \\
	   &&  -\left\langle x_{k+1}-x_k+\theta(x_k-x^*), A^T(y_k+\tfrac{1}{\theta}(y_{k+1}-y_k))\right\rangle \nonumber.
	   \end{eqnarray}

Since $F$ and $G$ are in symmetric positions in \eqref{eq_method} and  \( G \) is \( \mu_G \)-strongly convex and \( L_G \)-smooth, we can therefore, by employing computations similar to those used for $F$, derive the following analogous results for $G$:
\begin{eqnarray}\label{eqH3}
	&& G(y_{k+1})-G(y^*)-(1-\theta)(G(y_{k})-G(y^*)) \nonumber \\
	  && \leq  \tfrac{\theta(\theta^2-\mu_Gs)}{2s}\|\bar{y}_k-y^*\|^2+\tfrac{(1- \theta)}{2s}\|v_k\|^2- \tfrac{1}{2s}\|v_{k+1}\|^2 \\
	   && +\left\langle A\left(x_k+\tfrac{1}{\theta}(x_{k+1}-x_k)\right), y_{k+1}-y_k+\theta(y_k-y^*)\right\rangle.\nonumber
\end{eqnarray}
From \eqref{ques}, we have
\begin{eqnarray*}
&& \mathcal{L}(x_{k+1},y^*) -  \mathcal{L}(x^*,y_{k+1})\nonumber\\
&&\quad - (1-\theta)( \mathcal{L}(x_{k},y^*) -  \mathcal{L}(x^*,y_{k}))\nonumber\\
 && =  F(x_{k+1})-F(x^*)-(1-\theta)(F(x_{k})-F(x^*)) \\
 &&\quad +G(y_{k+1})-G(y^*)-(1-\theta)(G(y_{k})-G(y^*)) \\
 &&\quad +\left\langle x_{k+1}-x_k+\theta(x_k-x^*), A^T y^*\right\rangle\\
 &&\qquad-\left\langle Ax^*, y_{k+1}-y_k+\theta(y_k-y^*)\right\rangle.
\end{eqnarray*}	
Combining this  with \eqref{eq_24_1} and \eqref{eqH3}, we derive
\begin{eqnarray*}
&& \mathcal{L}(x_{k+1},y^*) -  \mathcal{L}(x^*,y_{k+1})\nonumber\\
&&\quad \leq (1-\theta)( \mathcal{L}(x_{k},y^*) -  \mathcal{L}(x^*,y_{k}))\nonumber\\
		 &&\qquad +\tfrac{\theta(\theta^2-\mu_Fr)}{2r}\|\bar{x}_k-x^*\|^2+ \tfrac{\theta(\theta^2-\mu_Gs)}{2s}\|\bar{y}_k-y^*\|^2\nonumber\\
 &&\qquad +\tfrac{(1- \theta)}{2r}\|u_k\|^2- \tfrac{1}{2r}\|u_{k+1}\|^2\\
	 &&\qquad +\tfrac{(1- \theta)}{2s}\|v_k\|^2- \tfrac{1}{2s}\|v_{k+1}\|^2\nonumber.
\end{eqnarray*}	
Since $\theta = \min\left\{ \sqrt{\mu_F r}, \sqrt{\mu_G s } \right\}$, it follows that ${\theta^2} \leq \mu_F r$ and ${\theta^2}\leq \mu_G s$.
Together with the definition of $\mathcal{E}_{k}$ in \eqref{eq_ener_dis} establishes the required result.	
\hspace*{\fill} \textbf{$\square$}
\end{proof}

 Based on Lemma \ref{lemmaL}, we can easily get the following linear convergence rate of Algorithm \ref{al_grad}.

\begin{theorem}\label{th_main}
		 Let $\{(x_k, y_k)\}_{k\geq 1}$ be the sequence generated by Algorithm \ref{al_grad}. Then,
		 \begin{eqnarray*}
			&&	\mathcal{L}(x_k,y^*)-\mathcal{L}(x^*,y_k)\leq {\mathcal{E}_1}(1-\theta)^{k-1},\\
			&&	\|x_k-x^*\|^2+\|y_k-y^*\|^2\leq \frac{2{\mathcal{E}_1}}{\min\{\mu_F,\mu_G\}}(1-\theta)^{k-1},
			\end{eqnarray*}
		where ${\mathcal{E}_1}=\mathcal{L}(x_1,y^*)-\mathcal{L}(x^*,y_1)+\tfrac{1}{2r}\left\|x_1\right\|^2+\tfrac{1}{2s}\left\|y_1\right\|^2$.
		\end{theorem}
\begin{proof}
It follows from Lemma \ref{lemmaL} that for any $k\geq 1$,
\[\mathcal{E}_{k+1} \leq  (1- \theta)\mathcal{E}_{k} \leq {\mathcal{E}_1} (1-\theta)^{k}.\]	
This together with the definition of $\mathcal{E}_{k}$ in \eqref{eq_ener_dis} yields
\begin{equation}\label{eq_he23_1}
		\mathcal{L}(x_k,y^*)-\mathcal{L}(x^*,y_k) \leq \mathcal{E}_k\leq {\mathcal{E}_1}(1-\theta)^{k-1}.
\end{equation}
From \eqref{eq_lag}, we have $-A^T y^*=\nabla F(x^*) ,\  Ax^*=\nabla G(y^*)$. Since $F$ and $G$ are strongly convex, it follows that
\begin{eqnarray*}
	&& \mathcal{L}(x_k,y^*)-\mathcal{L}(x^*,y_k) \\
	&&\quad= F(x_k)-F(x^*)+\langle A^Ty^*,  x_k-x^*\rangle\\
	&&\qquad + (G(y_k)-G(y^*)- \langle Ax^*,y_k-y^*\rangle)\\
	&&\quad\geq \tfrac{\mu_F}{2}\|x_k-x^*\|^2+\tfrac{\mu_G}{2}\|y_k-y^*\|^2.
\end{eqnarray*} 
Together with \eqref{eq_he23_1} implies the second inequality.
\hspace*{\fill} \textbf{$\square$}
\end{proof}

\begin{remark}
From Theorem \ref{th_main}, we obtain the 
\[\bigO\left(\left(1-\min\left\{\sqrt{\mu_F r}, \sqrt{\mu_G s}\right\}\right)^k\right)\]
  convergence rate for both the primal-dual gap and the iterative gap. Once we take $r=1/L_F$ and $s=1/L_G$, it achieves the linear convergence rate \eqref{opt_conv}, and yields oracle complexity \eqref{opt_oracle}, which improves upon those in the existing literature \citep{Borodich2023,Kovalev2022,Thekumparampil2022,Xie2020,Wang2020}. Furthermore,  our algorithm reduces to the NAG-\texttt{SC} method for unconstrained strongly convex optimization problems. Consequently, the convergence properties of Algorithm \ref{al_grad} also align with those of NAG-\texttt{SC}.
\end{remark}

\section{Continuous-time primal-dual dynamic}\label{sec3}

Continuous-time dynamical systems have long served as powerful theoretical frameworks for understanding and designing optimization algorithms. They provide continuous analogues of iterative schemes, enabling one to uncover the underlying convergence mechanisms through Lyapunov or energy-based analysis. In particular, the well-known NAG-\texttt{SC}~\eqref{al_Nag_sc} can be viewed as a time discretization of the second-order dynamical system~\eqref{dy_low_hb}, which explains its accelerated behavior and has inspired the development of various discrete accelerated methods. Recently, several studies have further investigated continuous-time dynamical systems and their accelerated counterparts for optimization problems with linear equality constraints~\citep{BotMP,BotJDE2022,HeSICON,LuoJota}, providing deeper insights into how damping  and rescaling contribute to the acceleration mechanism in both continuous and discrete settings.
For saddle-point problem~\eqref{ques}, a natural dynamic approach is the first-order primal-dual gradient flow~\citep{CherukuriSICON}:
\begin{equation}\label{dy_fd}
	\begin{cases}
		\dot{x}(t) + \nabla_x \mathcal{L}(x(t), y(t)) = 0, \\
		\dot{y}(t) - \nabla_y \mathcal{L}(x(t), y(t)) = 0,
	\end{cases}
\end{equation}
where the primal variable $x$ follows a gradient-descent direction and the dual variable $y$ follows gradient-ascent dynamics. While this dynamic coincides with the classical gradient flow for unconstrained problems, it does not exhibit accelerated convergence. Building on this, recent works~\citep{ZengIFAC,HeAMO,DingCoa} have proposed second-order primal-dual dynamical systems with different damping coefficients in the convex setting, achieving the optimal convergence rate of order~$\mathcal{O}(1/t^2)$. 
In this section, we make the following strongly convex assumptions about $F$ and $G$, which do not assume the global smoothness of $F$ and $G$.

\begin{assumption}\label{ass2}
$F:\mathbb{R}^n\to \mathbb{R}$ is $\mu_F$-strongly convex, $G:\mathbb{R}^m\to \mathbb{R}$ is $\mu_G$-strongly convex. The gradients $\nabla F$ and $\nabla G$ are locally Lipschitz continuous.
\end{assumption}

By employing Nesterov accelerated damping like in the dynamic \eqref{dy_low_hb},  we propose the following continuous-time primal-dual dynamic for solving problem \eqref{ques}:
\begin{equation*}
	\begin{cases}
		\ddot{x}(t)+2\sqrt{\mu_F}\dot{x}(t) +\nabla_x \mathcal{L}(x(t),y(t)+\gamma\dot{y}(t)) = 0, \\
		\ddot{y}(t)+2\sqrt{\mu_G}\dot{y}(t) -\nabla_y \mathcal{L}(x(t)+\gamma\dot{x}(t),y(t)) =0,
	\end{cases}
\end{equation*}
where $t\geq t_0\geq 0$ with initial conditions 
\[(x(t_0),y(t_0),\dot{x}(t_0),\dot{y}(t_0))=(x_0,y_0,u_0,v_0).\] 
Here, the extrapolation parameter  $\gamma$  is chosen as $\gamma = \max\{\tfrac{1}{\sqrt{\mu_F}},\tfrac{1}{\sqrt{\mu_G}}\}$.
Equivalently, the dynamical system can be expressed as
\begin{equation}\label{dynamic}
	\begin{cases}
		\ddot{x}(t)+2\sqrt{\mu_F}\dot{x}(t)  + \nabla F(x(t))+A^T(y(t)+\gamma\dot{y}(t)) = 0, \\
		\ddot{y}(t)+2\sqrt{\mu_G}\dot{y}(t) +\nabla G(y(t))-A(x(t)+\gamma\dot{x}(t))=0.
	\end{cases}
\end{equation}
  It  can provide theoretical insight into the acceleration mechanism underlying Algorithm~\ref{al_grad}.

From Assumption \ref{ass2}, we know that $\nabla F$ and $\nabla G$ are locally Lipschitz continuous, then by the Picard-Lindel{\"o}f theorem \citep[Theorem 2.2]{Teschl2012}, one can easily  get that for any  initial condition, there exists a unique solution $x(t) \in C^2([t_0, T), \mathbb{R}^n)$  and $y(t) \in C^2([t_0, T), \mathbb{R}^m)$ defined on a maximal interval $[t_0, T)$. If we can prove that $\sup_{t \in [t_0,T)} \|\dot{x}(t)\| < \infty$,
then, by the argument in \cite{HeTAC}, one can conclude that $T = +\infty$, which  prove the existence and uniqueness of a global solution. In the following analysis, we will verify that the $\sup_{t \in [t_0,T)} \|\dot{x}(t)\| < \infty$ holds. Consequently, throughout this section, we always assume that dynamic \eqref{dynamic} admits a global solution.

 To further investigate the long-time behavior of the trajectories, we construct a suitable energy function. Let $(x^*, y^*) \in \mathbb{R}^n \times \mathbb{R}^m$ denote the unique solution of  problem \eqref{ques}, and let $(x(t), y(t))$ be a global solution to the dynamic \eqref{dynamic}. Then, for all $t \geq t_0$, the saddle structure of the problem implies  $\mathcal{L}(x(t),y^*)-\mathcal{L}(x^*,y(t))\geq 0$. Next, we define the energy function  $\mathcal{E}:[t_0,+\infty)\to[0,+\infty)$ as
\begin{equation}\label{eq_energy}
	\mathcal{E}(t) =\gamma^2 (\mathcal{L}(x(t),y^*)-\mathcal{L}(x^*,y(t))) + \tfrac{1}{2}\|u(t)\|^2+\tfrac{1}{2}\|v(t)\|^2
\end{equation}
with
\begin{equation}\label{eq_ut}
	u(t) = x(t)-x^*+\gamma\dot{x}(t),
\end{equation}
and 
\begin{equation}\label{eq_vt}
	 v(t) = y(t)-y^*+\gamma\dot{y}(t).
\end{equation}
Compared with the discrete-time analysis in Section~\ref{sec2}, which focuses on establishing the  inequality
$\mathcal{E}_{k+1}\leq (1-\theta)\mathcal{E}_k$
for the energy sequence, the continuous-time analysis for dynamic~\eqref{dynamic} relies on constructing a differentiable Lyapunov function $\mathcal{E}(t)$~\eqref{eq_energy} and showing that its time derivative satisfies a differential inequality of the form
$\dot{\mathcal{E}}(t)\leq-{\mathcal{E}(t)}/{\gamma}$,
which directly implies exponential decay $\mathcal{E}(t)\leq \mathcal{E}(t_0)e^{-t/\gamma}$.
Hence, although both approaches aim to demonstrate linear convergence, the continuous-time proof is based on differential inequalities, whereas the discrete-time proof relies on recursive relations.

Next, we first show that the energy function $\mathcal{E}(t)$ is a nonincreasing function  under Assumption \ref{ass2}.

\begin{lemma}\label{lemma1}
Suppose that Assumption~\ref{ass2} holds and the energy function~$\mathcal{E}(t)$ is defined in~\eqref{eq_energy}. Then, for any $t \ge t_0$, 
\[\dot{\mathcal{E}}(t)\leq  -\tfrac{1}{\gamma}{\mathcal{E}}(t) -\tfrac{\gamma}{2}(\|\dot{x}(t)\|^2+\|\dot{y}(t)\|^2).\]
\end{lemma}
\begin{proof}
To differentiate $\tfrac{1}{2}\|u(t)\|^2$, we use the dynamic given by \eqref{dynamic} and \eqref{eq_ut}. This yields the following:
\begin{eqnarray}\label{eq_compu1}
&&	\tfrac{d}{dt}\left(\tfrac{1}{2}\|u(t)\|^2\right) = \langle u(t),\dot{u}(t)\rangle\nonumber \\
&& =  \langle u(t), \gamma(\ddot{x}(t)+2\sqrt{\mu_F}\dot{x}(t))+(1-2\gamma\sqrt{\mu_F})\dot{x}(t)\rangle\nonumber
\\
&& = -\gamma\langle u(t), \nabla F(x(t))+A^Ty^*+A^T(y(t)-y^*+\gamma\dot{y}(t))\rangle\nonumber\\
	&& \quad+(1-2\gamma\sqrt{\mu_F})\langle u(t), \dot{x}(t)\rangle.	  \end{eqnarray}
	  Next, we compute the right-hand side of the above equation. Since  $F$ is $\mu_F$-strongly convex, by the definition of $u(t)$ and $v(t)$, we have
\begin{eqnarray}\label{eq_compu2}
&&	\langle u(t),  \nabla F(x(t))+A^Ty^*\rangle\nonumber \\
&&\qquad =  \langle x(t)-x^*, \nabla F(x(t))+A^Ty^*\rangle \nonumber\\
	&&\qquad\quad +\gamma \langle\dot{x}(t), \nabla F(x(t))+A^Ty^*\rangle\\
	&&\qquad \geq F(x(t))-F(x^*)+\langle x(t)-x^*,A^Ty^*\rangle\nonumber\\
	&&\qquad\quad+\tfrac{\mu_F}{2}\|x(t)-x^*\|^2+\gamma \langle\dot{x}(t), \nabla F(x(t))+A^Ty^*\rangle\nonumber
\end{eqnarray}
and 
\begin{equation}\label{eq_compu3}
	\langle u(t), A^T(y(t)-y^*+\gamma\dot{y}(t))\rangle = \langle u(t),A^Tv(t)\rangle.
\end{equation}
Since $\gamma  = \max\left\{\tfrac{1}{\sqrt{\mu_F}},\tfrac{1}{\sqrt{\mu_G}}\right\}$, we have $1-2\gamma\sqrt{\mu_F}\leq -1$, and thus
\begin{eqnarray}\label{eq_compu4}
&& (1-2\gamma\sqrt{\mu_F})\langle u(t), \dot{x}(t)\rangle = \tfrac{1-2\gamma\sqrt{\mu_F}}{\gamma}\times \gamma\langle u(t), \dot{x}(t)\rangle\nonumber \\
&&\quad =  \tfrac{1-2\gamma\sqrt{\mu_F}}{\gamma}(\langle x(t)-x^*,  \gamma\dot{x}(t)\rangle + \gamma^2\|\dot{x}(t)\|^2)\nonumber\\
&&\quad =  \tfrac{1-2\gamma\sqrt{\mu_F}}{2\gamma}\|x(t)-x^*+\gamma\dot{x}(t)\|^2\ \\
&&\qquad +  \tfrac{1-2\gamma\sqrt{\mu_F}}{\gamma}\left(\tfrac{\gamma^2}{2}\|\dot{x}(t)\|^2-\tfrac{1}{2}\|x(t)-x^*\|^2\right)\nonumber \\
&&\quad \leq   -\tfrac{1}{2\gamma}\|u(t)\|^2-\tfrac{\gamma}{2}\|\dot{x}(t)\|^2-\tfrac{1-2\gamma\sqrt{\mu_F}}{2\gamma}\|x(t)-x^*\|^2,\nonumber
\end{eqnarray}
where the third equality follows from  the identity:
\[ \langle x,y\rangle = \tfrac{1}{2}(\|x+y\|^2-\|x\|^2-\|y\|^2), \quad \forall x,y\in\mathbb{R}^n.\]
Combining \eqref{eq_compu1}-\eqref{eq_compu4} with  $\mu_F\gamma^2-2\sqrt{\mu_F}\gamma+1\geq 0$ for any $\gamma$, we obtain
\begin{eqnarray}\label{eq_compu}
&&	\tfrac{d}{dt}\left(\tfrac{1}{2}\|u(t)\|^2\right) \leq   -\tfrac{1}{2\gamma}\|u(t)\|^2-\tfrac{\gamma}{2}\|\dot{x}(t)\|^2\nonumber\\ 
&&\quad	-\gamma(F(x(t))-F(x^*)+\langle x(t)-x^*,A^Ty^*\rangle)\nonumber\\
		&&\quad-\gamma\langle u(t), A^Tv(t)\rangle -\gamma^2 \langle\dot{x}(t), \nabla F(x(t))+A^Ty^*\rangle.
	\end{eqnarray}

Now, we differentiate  $\tfrac{1}{2}\|v(t)\|^2$. 
By similar computations as for \eqref{eq_compu1}-\eqref{eq_compu4}, we obtain
\begin{eqnarray}\label{eq_compv}
		&&\tfrac{d}{dt}\left(\tfrac{1}{2}\|v(t)\|^2\right)\leq -\tfrac{1}{2\gamma}\|v(t)\|^2-\tfrac{\gamma}{2}\|\dot{y}(t)\|^2\nonumber\\
		&&\quad -\gamma(G(y(t))-G(y^*)-\langle Ax^*,y(t)-y^*\rangle)\nonumber\\
		&&\quad +\gamma\langle Au(t), v(t)\rangle-\gamma^2 \langle\dot{y}(t), \nabla G(y(t))-Ax^*\rangle.
\end{eqnarray}

Combining \eqref{eq_energy}, \eqref{eq_compu}, and \eqref{eq_compv}, we compute
\begin{eqnarray*}
	\dot{\mathcal{E}}(t) &=& \gamma^2\langle \dot{x}(t), \nabla F(x(t))+A^Ty^*\rangle \nonumber\\
	&&-\gamma^2\langle \dot{y}(t), -\nabla G(x(t))+Ax^*\rangle\nonumber \\
	&&+\tfrac{d}{dt}\left(\tfrac{1}{2}\|u(t)\|^2\right)+\tfrac{d}{dt}\left(\tfrac{1}{2}\|v(t)\|^2\right)\nonumber \\
&\leq& -\gamma(F(x(t))+G(y(t))-F(x^*)-G(y^*)\nonumber\\
&&+\langle x(t),A^Ty^*\rangle - \langle Ax^*,y(t)\rangle)\\
&& -\tfrac{1}{2\gamma}(\|u(t)\|^2+\|v(t)\|^2 )-\tfrac{\gamma}{2}(\|\dot{x}(t)\|^2+\|\dot{y}(t)\|^2).\nonumber
\end{eqnarray*}
This, together with \eqref{eq_energy} yields the result. \hspace*{\fill} \textbf{$\square$}
\end{proof}

Now, we are in a position to investigate the convergence properties of the trajectories.

\begin{theorem}\label{th_res}
	Suppose that Assumption \ref{ass2} holds. Then the trajectory $(x(t),y(t))$ of the dynamic \eqref{dynamic} satisfies the following convergence results:
	\begin{itemize}
	 \setlength{\itemsep}{3pt}
	\item [(i)] 
	$\mathcal{L}(x(t),y^*)-\mathcal{L}(x^*,y(t))\leq \mathcal{O}(e^{-\min\{\sqrt{\mu_F},\sqrt{\mu_G}\}t})$.\
	 \item [(ii)] $\|x(t)-x^*\|^2+\|y(t)-y^*\|^2\leq \mathcal{O}(e^{-\min\{\sqrt{\mu_F},\sqrt{\mu_G}\}t}).$
	   \item [(iii)] $\|\dot{x}(t)\|^2+\|\dot{y}(t)\|^2\leq \mathcal{O}(e^{-\min\{\sqrt{\mu_F},\sqrt{\mu_G}\}t})$.
\end{itemize}
\end{theorem}
\begin{proof}
From Lemma \ref{lemma1}, we know that
\[\dot{\mathcal{E}}(t)\leq  -\tfrac{1}{\gamma}{\mathcal{E}}(t) -\tfrac{\gamma}{2}(\|\dot{x}(t)\|^2+\|\dot{y}(t)\|^2)\]
for any $t\in[t_0,+\infty)$. Multiply both sides of this inequality by  $e^{t/\gamma}$, we have 
\[
	 \tfrac{d}{dt}\left(e^{t/\gamma} {\mathcal{E}}(t)\right)\leq - \tfrac{\gamma e^{t/\gamma}}{2}(\|\dot{x}(t)\|^2+\|\dot{y}(t)\|^2)\leq 0.
\]
This implies that ${\mathcal{E}}(t) \leq  \tfrac{e^{t_0/\gamma} {\mathcal{E}}(t_0)}{ e^{t/\gamma}}$ for any $t\geq t_0$. 
By the definition of ${\mathcal{E}}(t)$, we have 
\[
	\mathcal{L}(x(t),y^*)-\mathcal{L}(x^*,y(t)) 
	\leq  \tfrac{e^{t_0/\gamma} {\mathcal{E}}(t_0)}{\gamma^2 e^{t/\gamma}}  =\mathcal{O}(e^{-t/\gamma})
\]
and
\[\|u(t)\|^2+\|v(t)\|^2\leq  \tfrac{2e^{t_0/\gamma} {\mathcal{E}}(t_0)}{ e^{t/\gamma}}  =\mathcal{O}(e^{-t/\gamma}). \]
Together with strongly convex assumption of $F$ and $G$  implies
\[\|x(t)-x^*\|^2\leq \mathcal{O}(e^{-t/\gamma}),\quad  \|y(t)-y^*\|^2\leq \mathcal{O}(e^{-t/\gamma}), \]
and then
\[\|\dot{x}(t)\|^2\leq \tfrac{2}{\gamma^2}(\|u(t)\|^2+\|x(t)-x^*\|^2) =\mathcal{O}(e^{-t/\gamma}), \]
\[\|\dot{y}(t)\|^2\leq \tfrac{2}{\gamma^2}(\|v(t)\|^2+\|y(t)-y^*\|^2) =\mathcal{O}(e^{-t/\gamma}).\]
Combining these results with $\gamma = \max\left\{{1}/{\sqrt{\mu_F}},{1}/{\sqrt{\mu_G}}\right\}$ yields $(i)-(iii)$.\hspace*{\fill} \textbf{$\square$}
\end{proof}

\begin{remark}
Under Assumption \ref{ass2}, there must exist a unique local solution $x(t) \in C^2([t_0,T),\mathbb{R}^n),\ y(t) \in C^2([t_0,T),\mathbb{R}^m)$, defined on a maximal interval $[t_0,T)$. Define the energy function $\mathcal{E}:[t_0,T)\to[0,+\infty)$ as given in \eqref{eq_energy}. By Theorem~\ref{th_res} $(iii)$, it follows that
$\sup_{t \in [t_0,T)} \|\dot{x}(t)\| < +\infty$. Invoking similar arguments as in \cite{HeTAC}, we  deduce that $T = +\infty$, which ensures the global existence of the solution.
\end{remark}

\begin{remark}
Under the strongly convex-strongly concave assumption (Assumption \ref{ass2}), we show that the accelerated primal-dual dynamic \eqref{dynamic} exhibits a $\mathcal{O}(e^{-\min\{\sqrt{\mu_F},\sqrt{\mu_G}\}t})$  linear convergence rate for the primal-dual gap.  If we set $A=0$, problem \eqref{ques} reduces to two separate unconstrained strongly convex optimization problems. By focusing solely on the $F$-part (in this case, we can choose $G$ to be any strongly convex function, such as $G(y)=\tfrac{\mu_F}{2}\|y\|^2)$, the obtained convergence rate of  $\mathcal{O}(e^{-\min\{\sqrt{\mu_F},\sqrt{\mu_G}\}t})$  simplifies to $\mathcal{O}(e^{-\sqrt{\mu_F}t})$, which recovers the convergence results from the Nesterov accelerated dynamical system \eqref{dy_low_hb} in \citet{LuoMP,Wilson,Siegel}.
 \end{remark}

Although Algorithm~\ref{al_grad} and the dynamical system~\eqref{dynamic} can both be regarded as extensions of the discrete- and continuous-time Nesterov accelerated gradient methods for unconstrained optimization problems, the relationship between them is not immediately clear. To gain a deeper understanding of their connection, we next provide a detailed derivation showing how the discrete Algorithm~\ref{al_grad} can be viewed as a time discretization of  the rescaled version of the continuous-time system~\eqref{dynamic}. To this end, consider the following rescaled version of the dynamical system~\eqref{dynamic}:
\begin{equation}\label{dy_2}
\begin{cases}
\ddot{x}(t)+2\sqrt{\alpha}\dot{x}(t) +\beta_1\nabla_x \mathcal{L}(x(t),y(t)+\gamma\dot{y}(t)) = 0, \\
\ddot{y}(t)+2\sqrt{\alpha}\dot{y}(t) -\beta_2\nabla_y \mathcal{L}(x(t)+\gamma\dot{x}(t),y(t)) = 0,
\end{cases}
\end{equation}
where $\sqrt{\alpha} = \tfrac{1}{\gamma} = \min\{\sqrt{\mu_F}, \sqrt{\mu_G}\}$, and $\beta_1, \beta_2 > 0$ are two scaling coefficients.
Under Assumption~\ref{ass2}, both $F$ and $G$ are $\alpha$-strongly convex.
We define the Lyapunov (energy) function $
	\mathcal{E}(t) =\gamma^2 (\mathcal{L}(x(t),y^*)-\mathcal{L}(x^*,y(t)))+ \tfrac{1}{2\beta_1}\|u(t)\|^2+\tfrac{1}{2\beta_2}\|v(t)\|^2$,
where $u(t)$ and $v(t)$ are given in~\eqref{eq_ut} and~\eqref{eq_vt}, respectively.
By following similar steps to those in the proof of Theorem~\ref{th_res}, it can be shown that the rescaled dynamical system~\eqref{dy_2} achieves the same exponential convergence rate as~\eqref{dynamic}, that is $\mathcal{O}(e^{-\min\{\sqrt{\mu_F}, \sqrt{\mu_G}\}t})$.

To further establish the link between the continuous and discrete dynamics, we now perform a time discretization of~\eqref{dy_2}.
Let the step size be $\sqrt{h}$, satisfying $\tfrac{\gamma}{\sqrt{h}} = \tfrac{1}{\theta}$, where $\theta$ denotes the parameter used in Algorithm~\ref{al_grad}. Since $\sqrt{\alpha} = \tfrac{1}{\gamma}$, it follows that $\sqrt{\alpha h} = \theta$. We define the discrete time instants $t_k = k\sqrt{h}$ and denote $x_k = x(t_k)$ and $y_k = y(t_k)$.
Applying an explicit discretization scheme to~\eqref{dy_2} with respect to $F$ and $G$ gives
\begin{equation*}
\begin{cases}{}
	\tfrac{x_{k+1}-2x_k+x_{k-1}}{h} +\sqrt{\alpha}\tfrac{x_{k+1}-x_k}{\sqrt{h}}+\sqrt{\alpha}\tfrac{x_{k}-x_{k-1}}{\sqrt{h}}&\\
	\quad+\beta_1(\nabla F(\bar{x}_k)+A^T(y_{k}+\tfrac{\gamma}{\sqrt{h}}(y_{k+1}-y_{k})) &= 0,\\
	\tfrac{y_{k+1}-2y_k+y_{k-1}}{h} +\sqrt{\alpha}\tfrac{y_{k+1}-y_k}{\sqrt{h}}+\sqrt{\alpha}\tfrac{y_{k}-y_{k-1}}{\sqrt{h}}&\\
	\quad+\beta_2(\nabla G(\bar{y}_k)-A(x_{k+1}+\tfrac{\gamma}{\sqrt{h}}(x_{k}-x_{k})) &= 0,
 \end{cases}
\end{equation*}
where $(\bar{x}_k, \bar{y}_k) = (x_k, y_k) + \tfrac{1 - \theta}{1 + \theta} \left[ (x_k, y_k) - (x_{k-1}, y_{k-1}) \right]$ is an extrapolated point used to introduce Nesterov-type acceleration. After rearranging terms and simplifying, the above equation leads to the following discrete  rule:
\begin{equation*}
\begin{cases}{}
	x_{k+1} = \bar{x}_k - \tfrac{\beta_1 h}{1+\theta} \left( \nabla F(\bar{x}_k) + A^T \left( y_k + \tfrac{1}{\theta}(y_{k+1} - y_k) \right) \right),  \\
	y_{k+1} = \bar{y}_k - \tfrac{\beta_2 h}{1+\theta} \left( \nabla G(\bar{y}_k) - A \left( x_k + \tfrac{1}{\theta}(x_{k+1} - x_k) \right) \right). \end{cases}
\end{equation*}
By setting  scaling coefficients as $\beta_1 = \tfrac{r(1 + \theta)}{h}$ and $\beta_2 = \tfrac{s(1 + \theta)}{h}$,  we recover the iteration scheme in \eqref{eq_method}, which corresponds exactly to the Algorithm \ref{al_grad}. This shows the connection between Algorithm \ref{al_grad}  and dynamic \eqref{dynamic}.
 
\begin{figure}[htbp]
\centering

\begin{minipage}{\linewidth}
\centering
\includegraphics[width=0.48\linewidth]{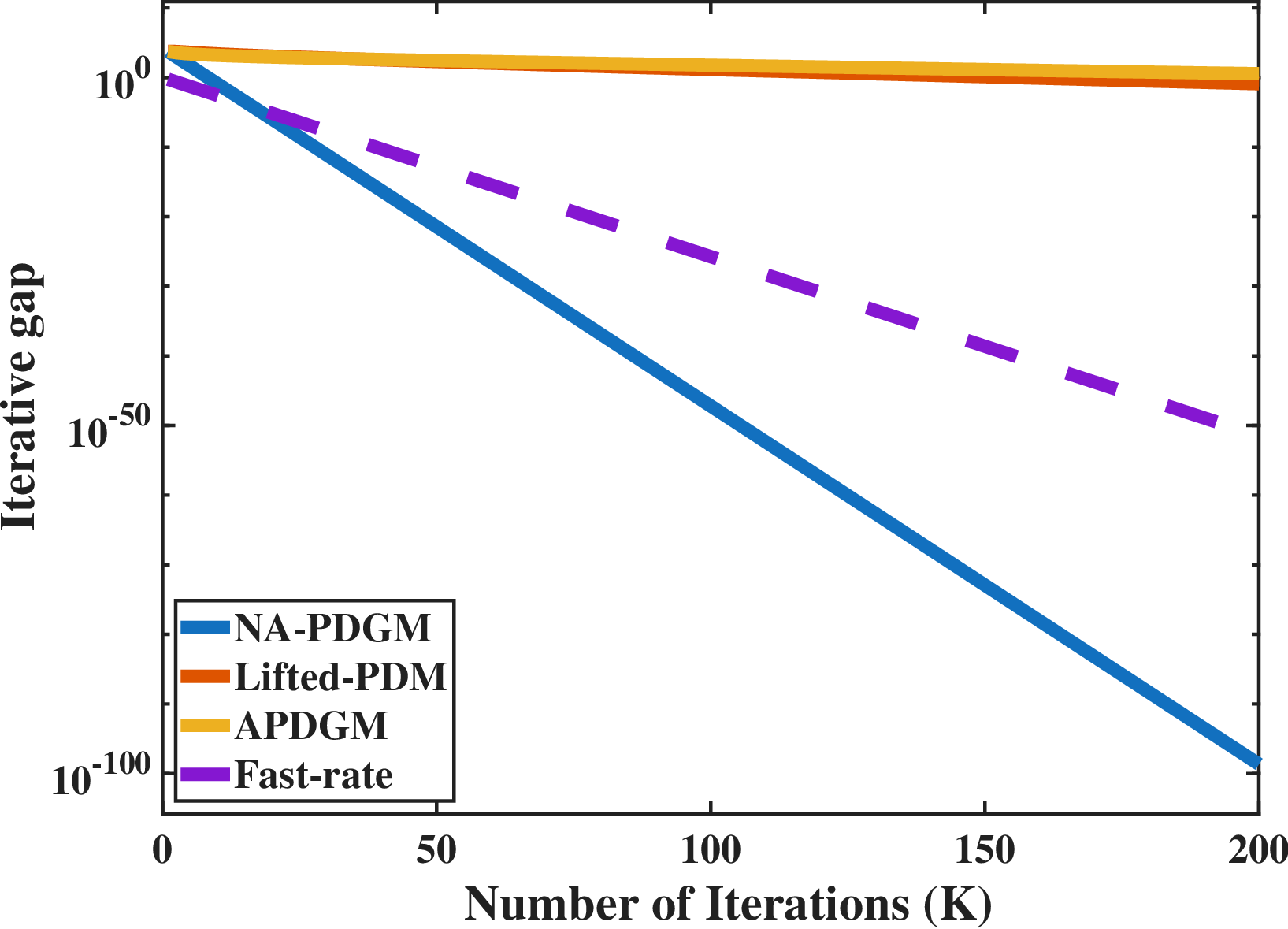}\hfill
\includegraphics[width=0.48\linewidth]{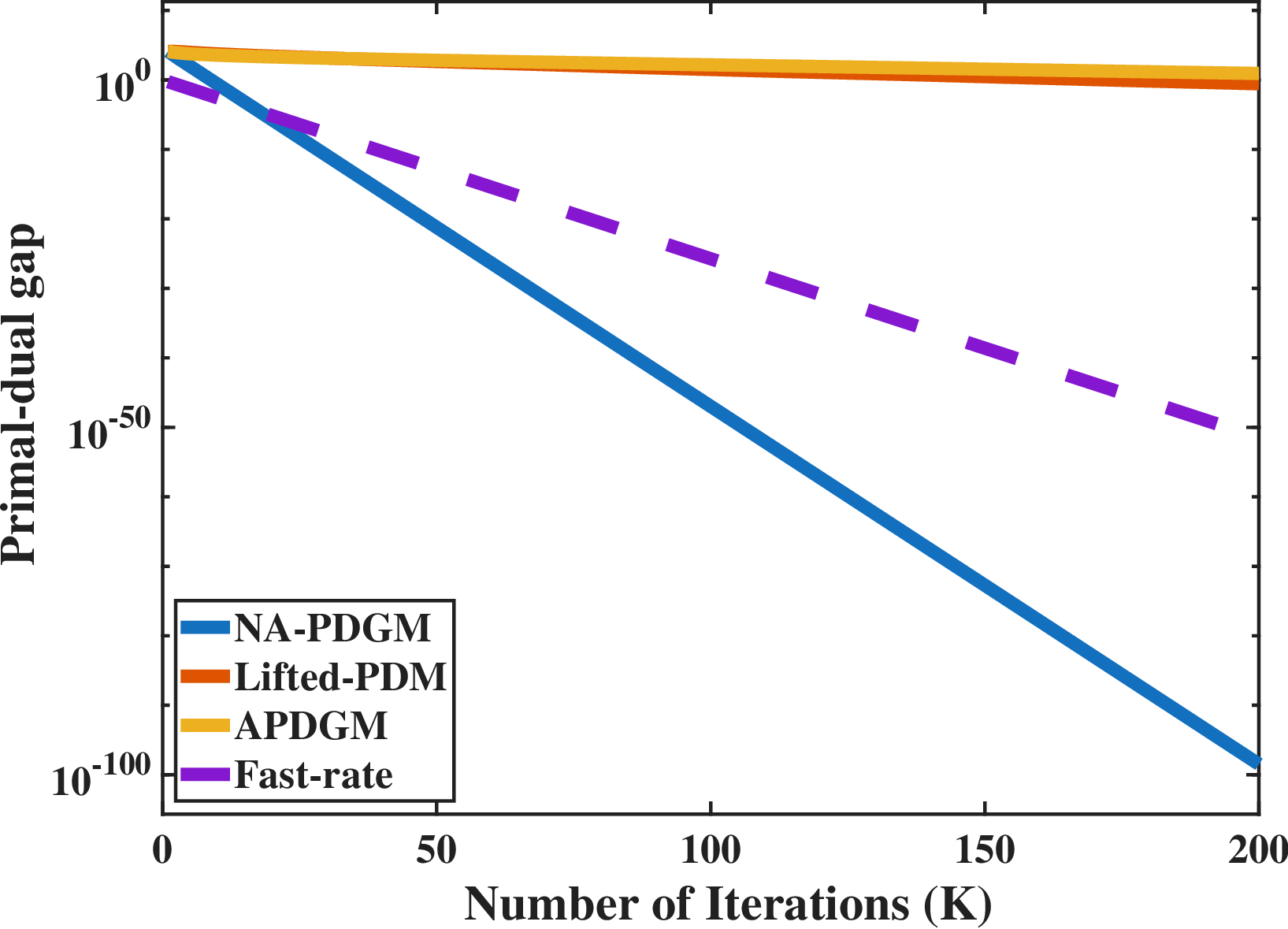}

\vspace{0.2cm}
\centerline{Case 1: $\kappa(R)=2$ and $\kappa(S)=5$ ($\theta=0.4472$)}
\end{minipage}

\vspace{0.35cm}
\begin{minipage}{\linewidth}
\centering
\includegraphics[width=0.48\linewidth]{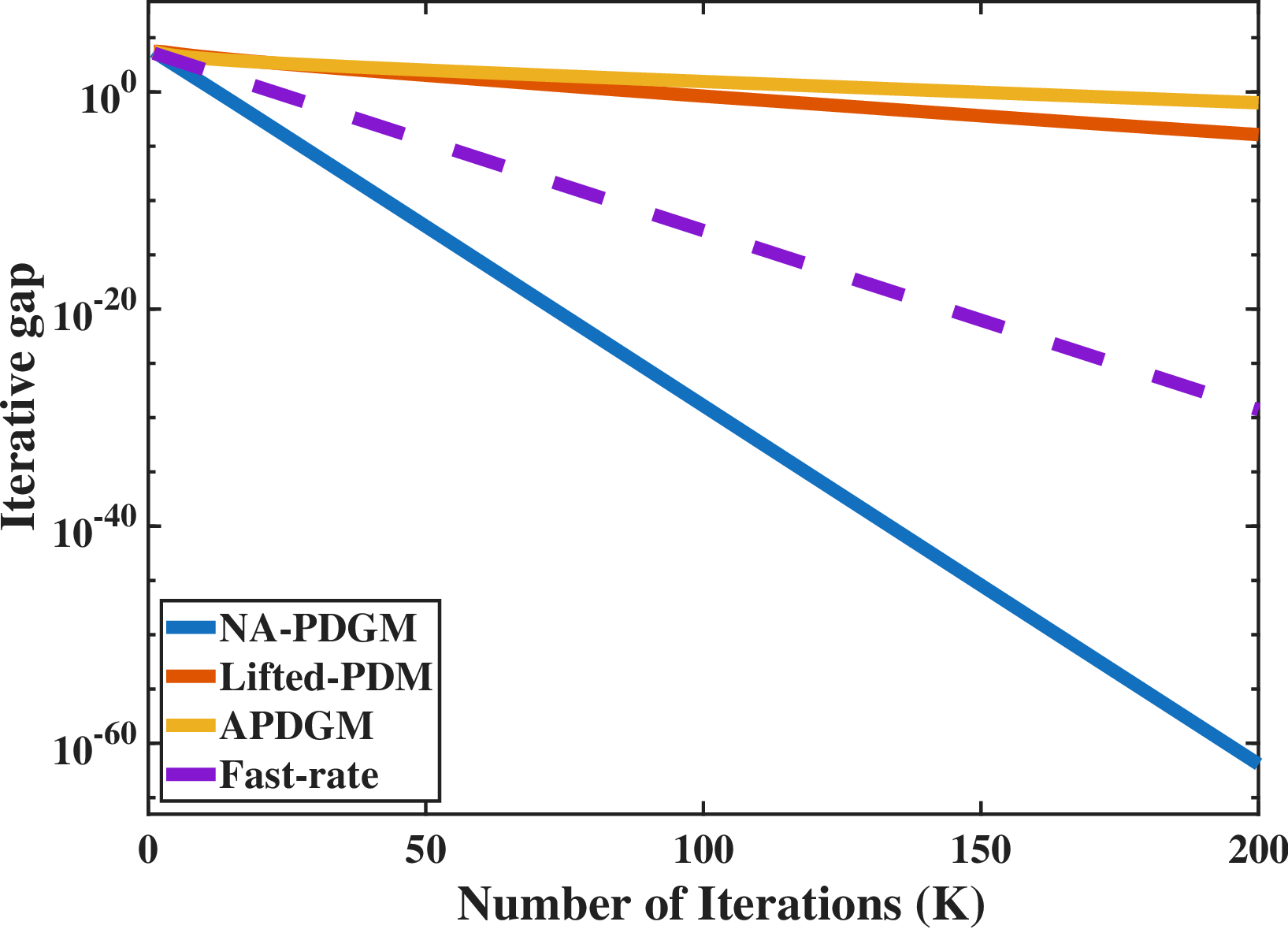}\hfill
\includegraphics[width=0.48\linewidth]{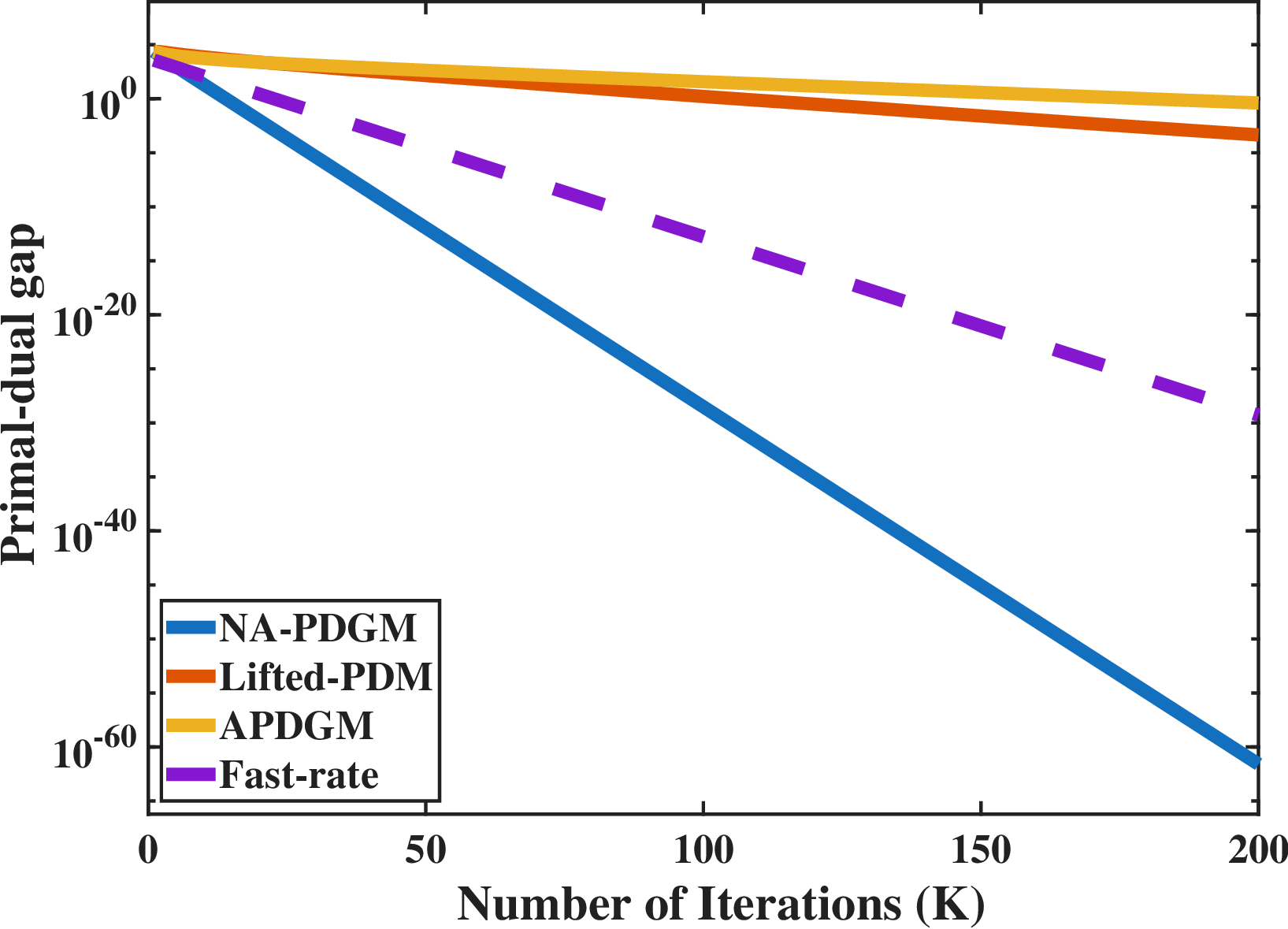}

\vspace{0.2cm}
\centerline{Case 1: $\kappa(R)=5$ and $\kappa(S)=10$ ($\theta=0.3162$)}
\end{minipage}

\vspace{0.35cm}

\begin{minipage}{\linewidth}
\centering
\includegraphics[width=0.48\linewidth]{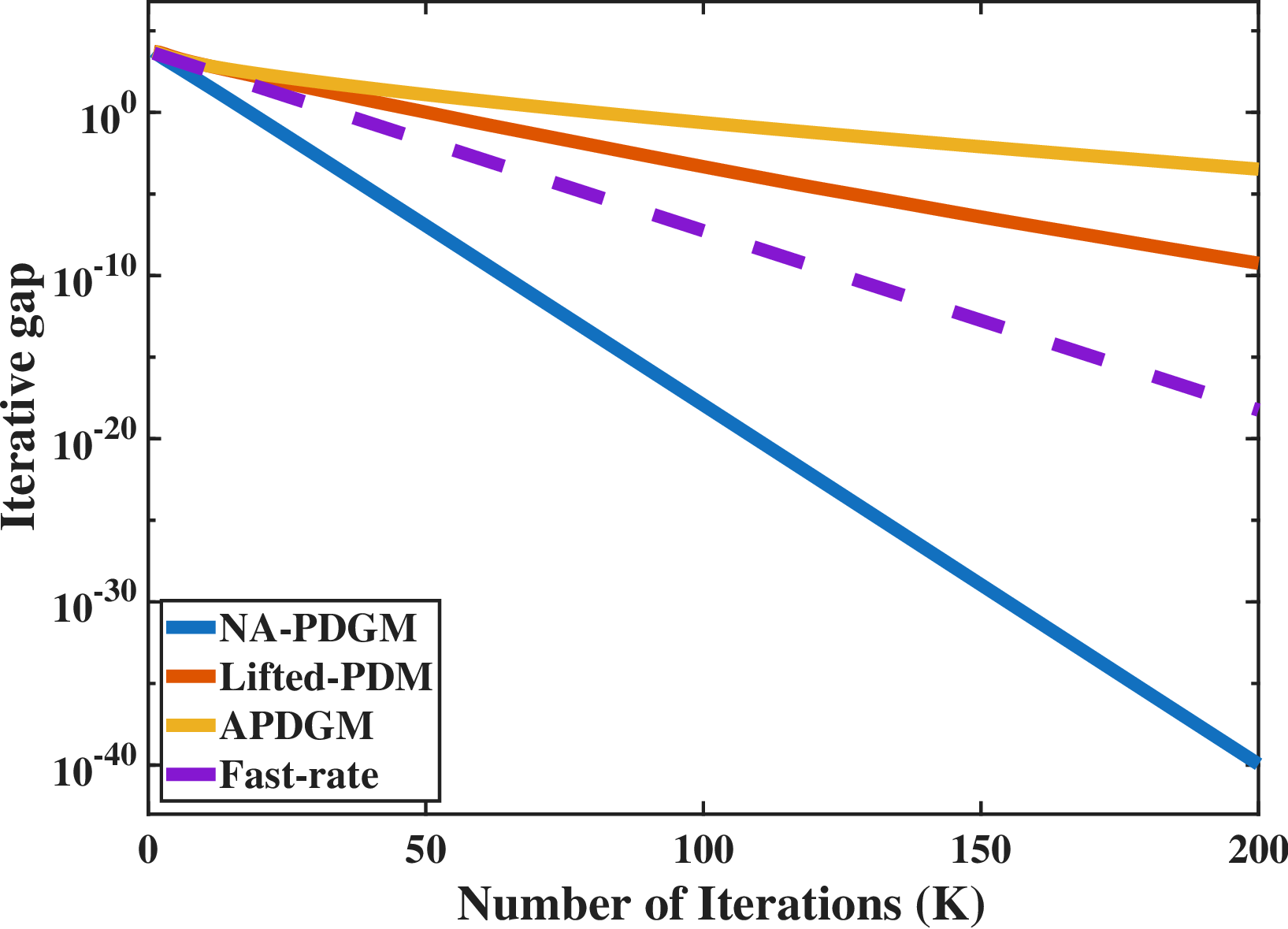}\hfill
\includegraphics[width=0.48\linewidth]{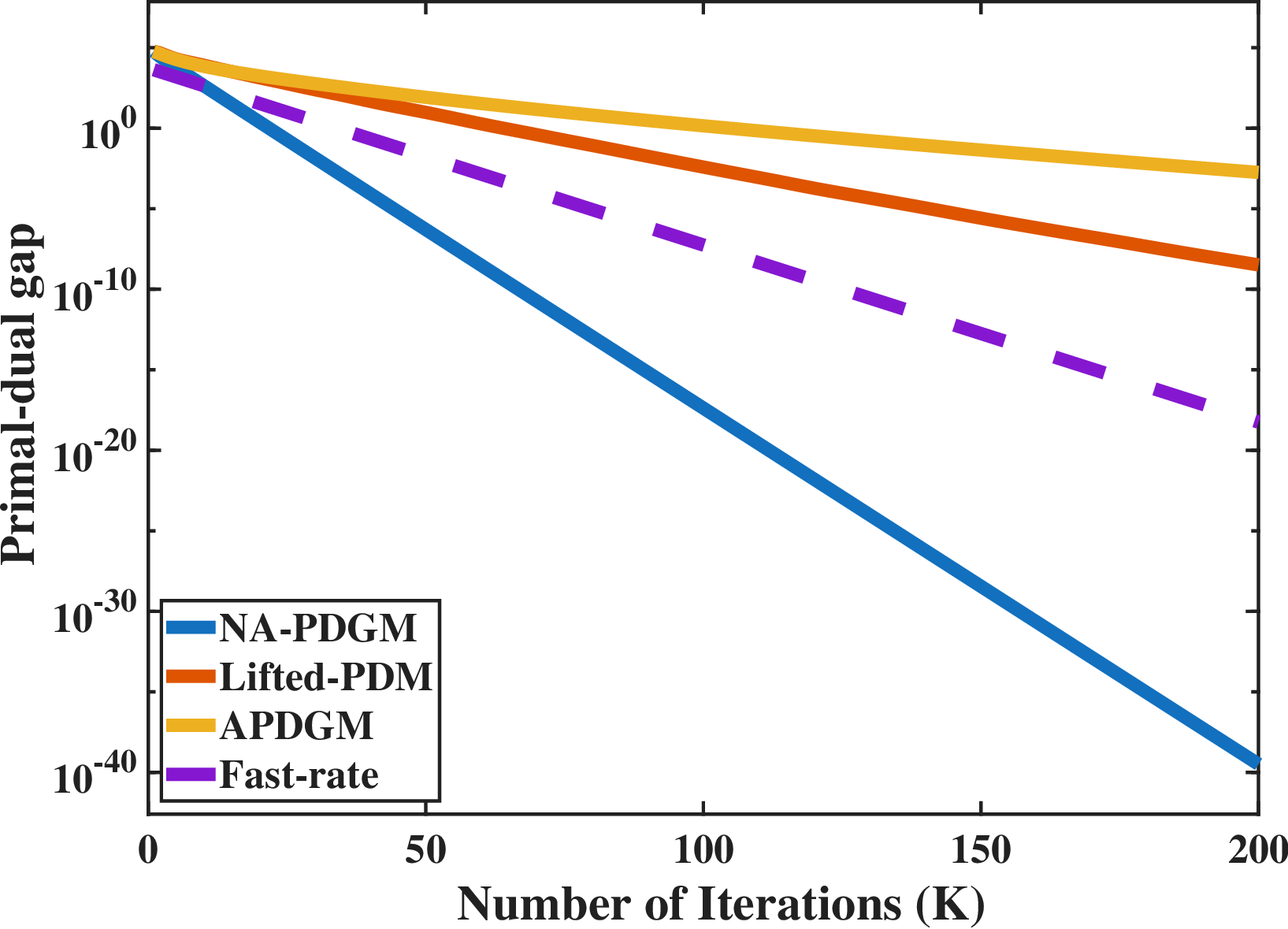}

\vspace{0.2cm}
\centerline{Case 2: $\kappa(R)=20$ and $\kappa(S)=20$ ($\theta=0.2236$)}
\end{minipage}
\caption{Numerical results of algorithms for Example 1}\label{fig1}
\end{figure}

\section{Numerical experiments}\label{sec4}

\emph{Example 1.} Consider the quadratic problem of the form \eqref{ques_quad}.  We randomly generate symmetric positive matrixes $R\in\mathbb{R}^{n\times n}$ and $S\in\mathbb{R}^{m\times m}$.  We generate the matrix $A\in\mathbb{R}^{m \times n}$ from a standard Gaussian distribution. In this case, we can easily that $(x^*,y^*)=(0,0)$ is the unique solution of problem \eqref{ques_quad}. In the numerical experiments, we set $m=3000$ and $n=2500$, and choose $\lambda_{\min}(R)=\lambda_{\min}(S)=1$. We consider different values of $\kappa(R)$ and $\kappa(S)$, which correspond to different choices of $\theta$ in Algorithm \ref{al_Nag_sc}. We compare the performance of the Nesterov Accelerated Primal-Dual Gradient Method (NA-PDGM, namely Algorithm \ref{al_grad} with $r=1/(2\lambda_{\max}(R))$ and $s=1/(2\lambda_{\max}(S))$), the Lifted Primal-Dual Method (Lifted-PDM) in \citet{Thekumparampil2022}, the Accelerated Primal-Dual Gradient Method (APDGM) in \citet{Kovalev2022}, and the fast rate $(1-\min\{1/\sqrt{\kappa(R)},\,1/\sqrt{\kappa(S)}\})^k$ given in Theorem \ref{th_main}. In Fig.~\ref{fig1}, we plot the primal-dual gap $\mathcal{L}(x_k,y^*)-\mathcal{L}(x^*,y_k)$ and the iterative gap $\|x_k-x^*\|^2+\|y_k-y^*\|^2$ versus the number of iterations $K$ for the different algorithms. The results show that NA-PDGM achieves faster linear convergence than the other methods under different condition numbers of $R$ and $S$. Moreover, the numerical results are consistent with the linear convergence rate predicted by Theorem \ref{th_main}, and among all the compared methods, only NA-PDGM attains this fast convergence rate.

\begin{figure}[htbp]
\centering

\begin{minipage}{\linewidth}
\centering
\includegraphics[width=0.48\linewidth]{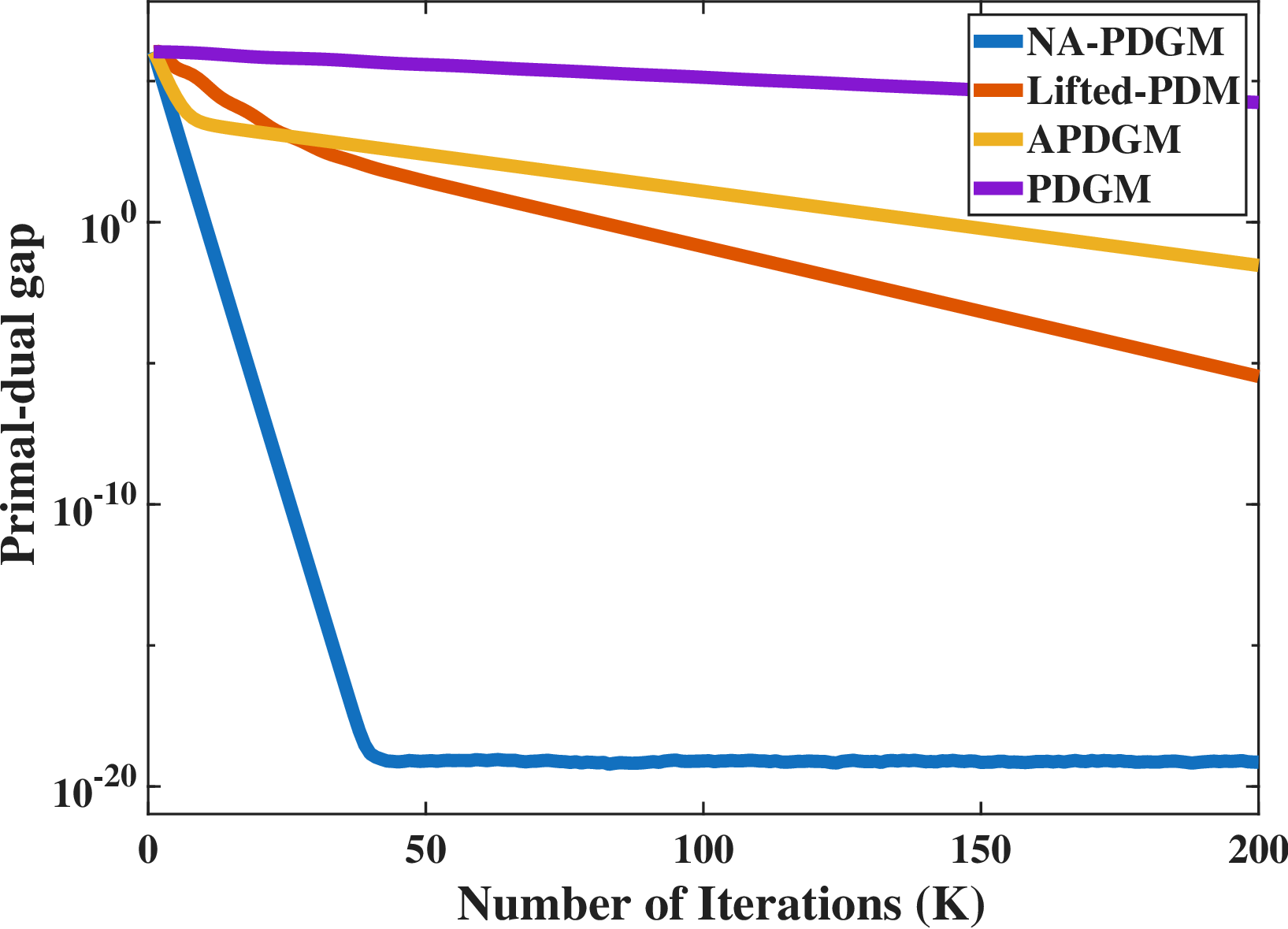}\hfill
\includegraphics[width=0.48\linewidth]{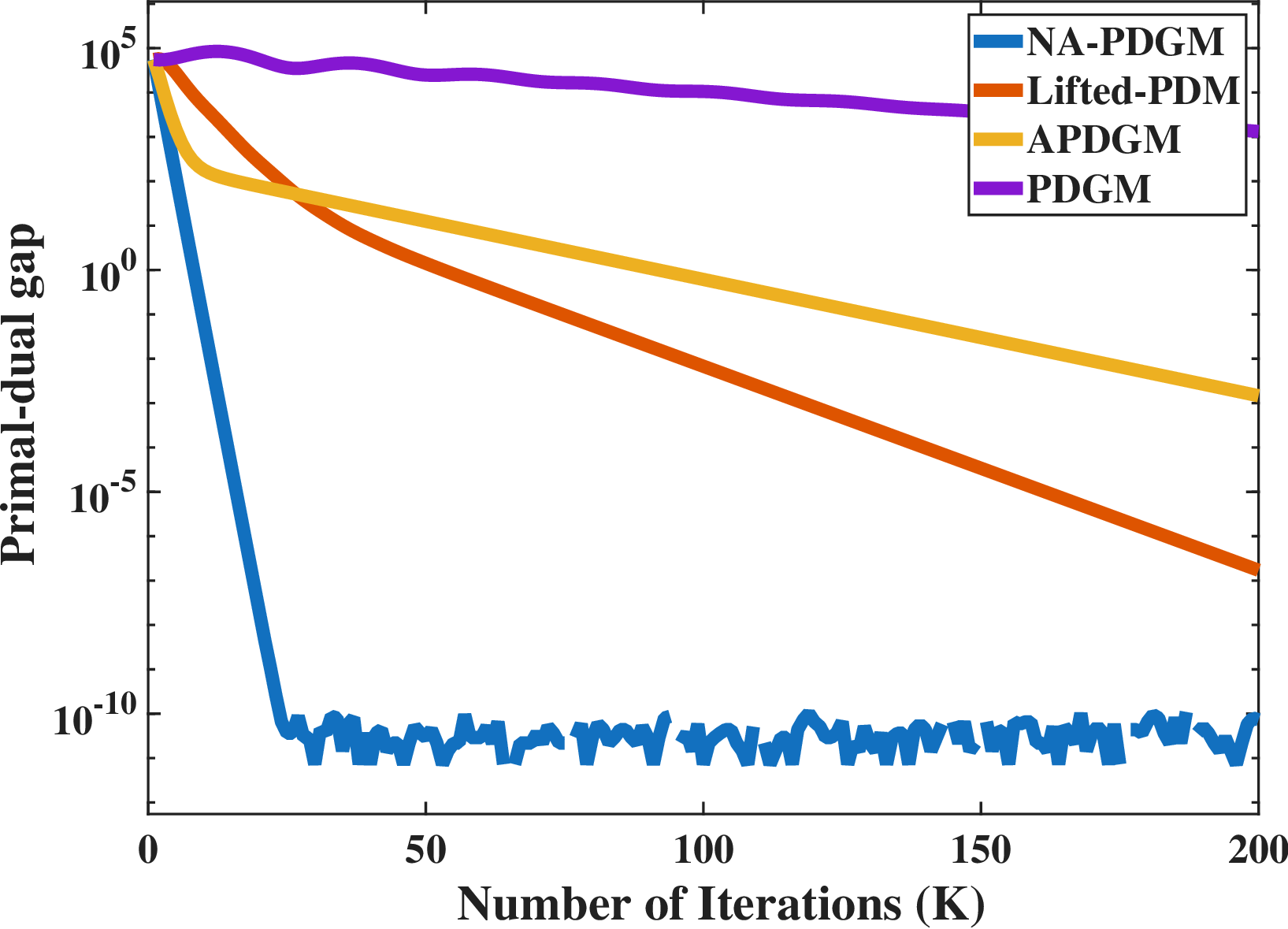}

\vspace{0.2cm}
\centerline{Case 1: $u=0.1$ and $r=0.2$ ($\theta=0.5345$)}
\end{minipage}

\vspace{0.35cm}
\begin{minipage}{\linewidth}
\centering
\includegraphics[width=0.48\linewidth]{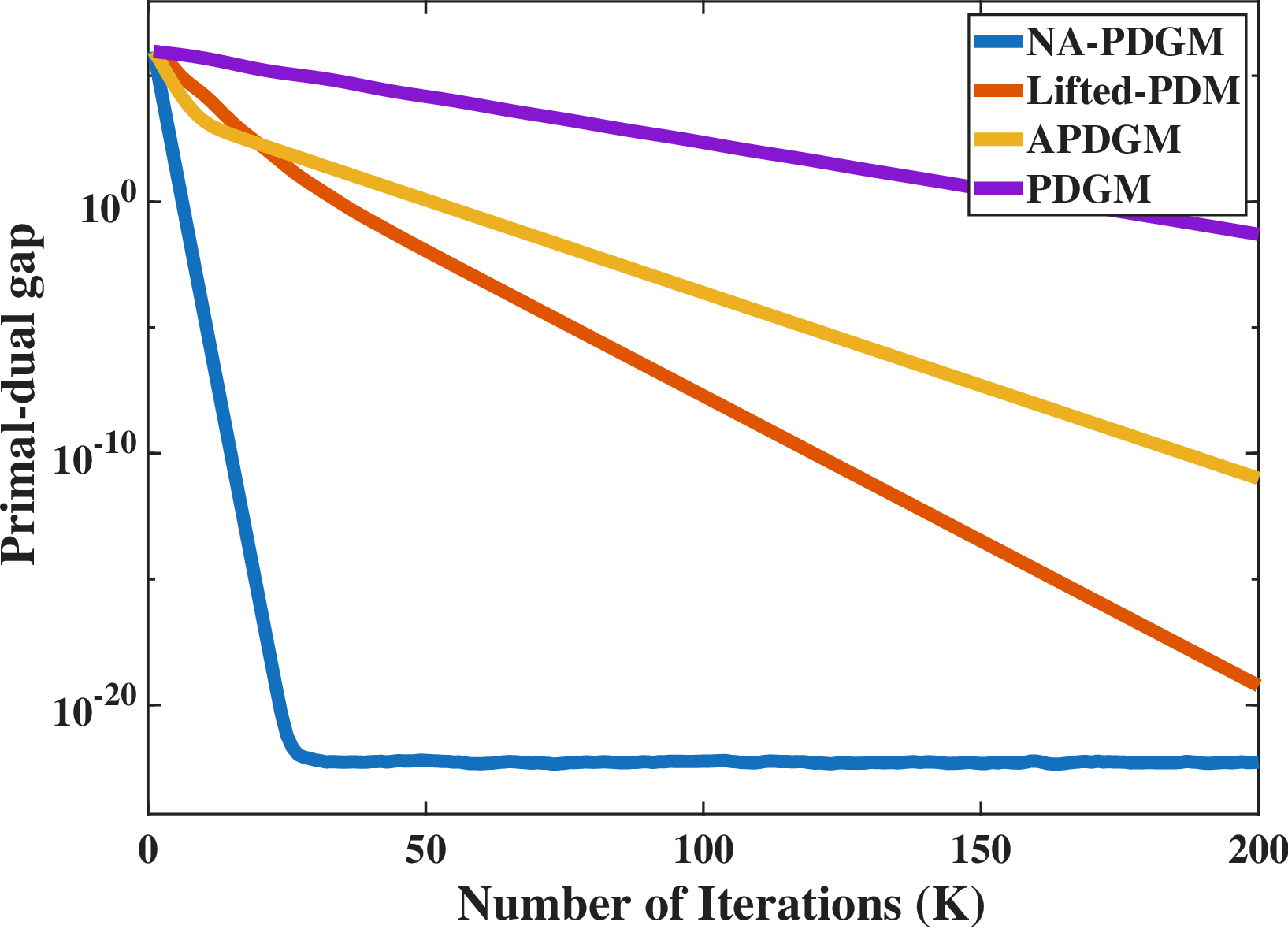}\hfill
\includegraphics[width=0.48\linewidth]{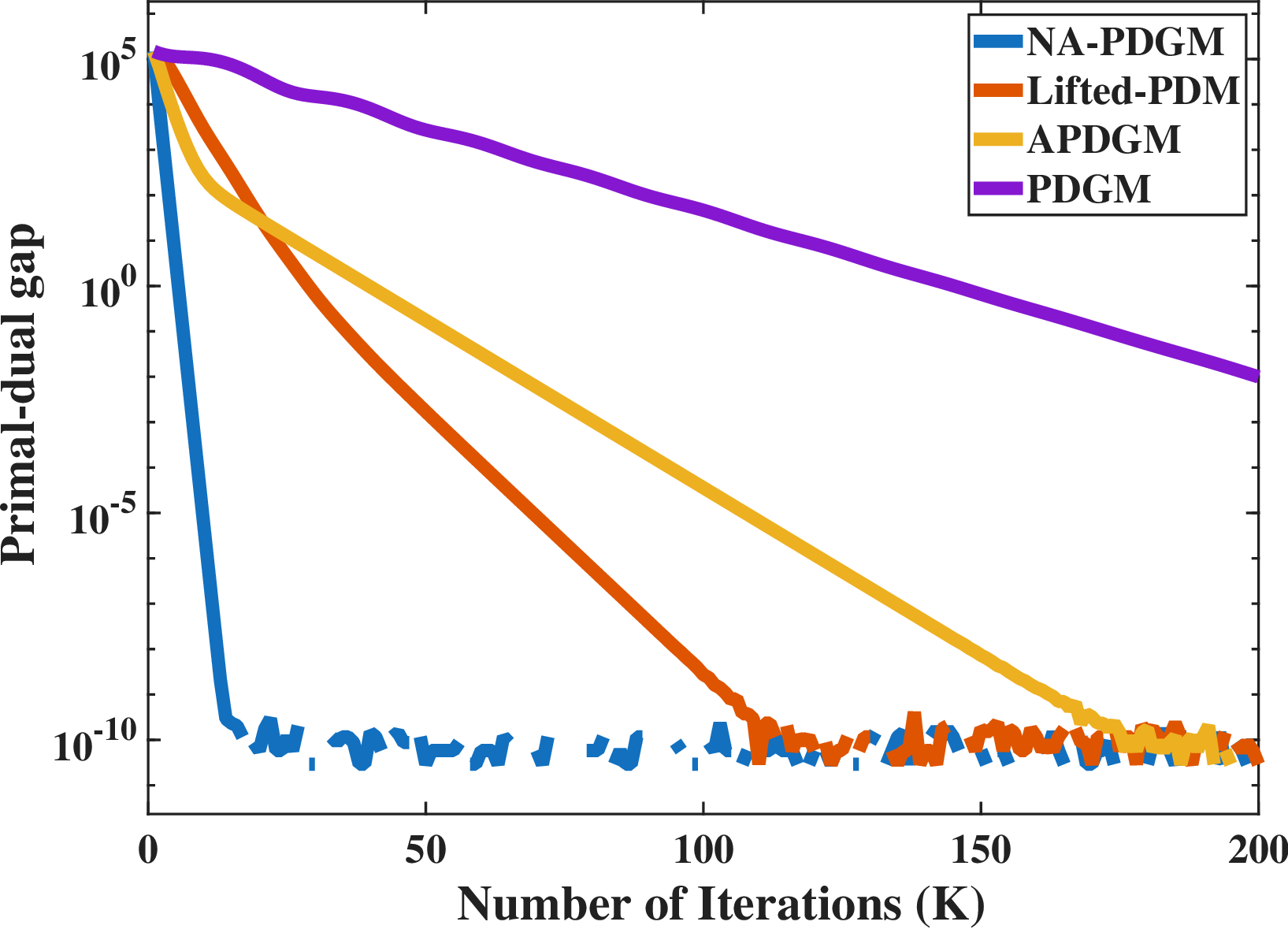}

\vspace{0.2cm}
\centerline{Case 3: $u=0.3$ and $r=0.5$ ($\theta=0.7385$)}
\end{minipage}

\vspace{0.35cm}

\begin{minipage}{\linewidth}
\centering
\includegraphics[width=0.48\linewidth]{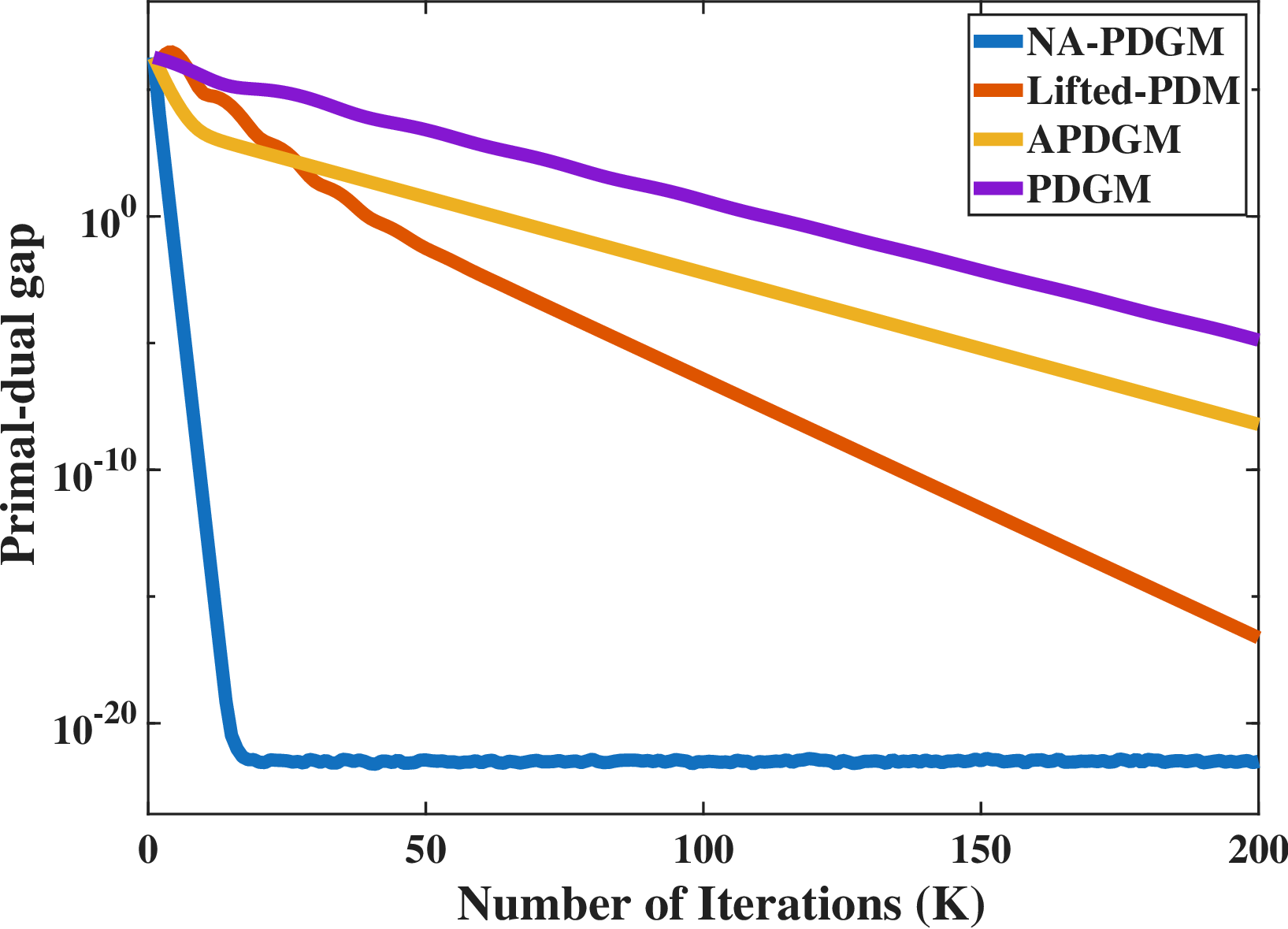}\hfill
\includegraphics[width=0.48\linewidth]{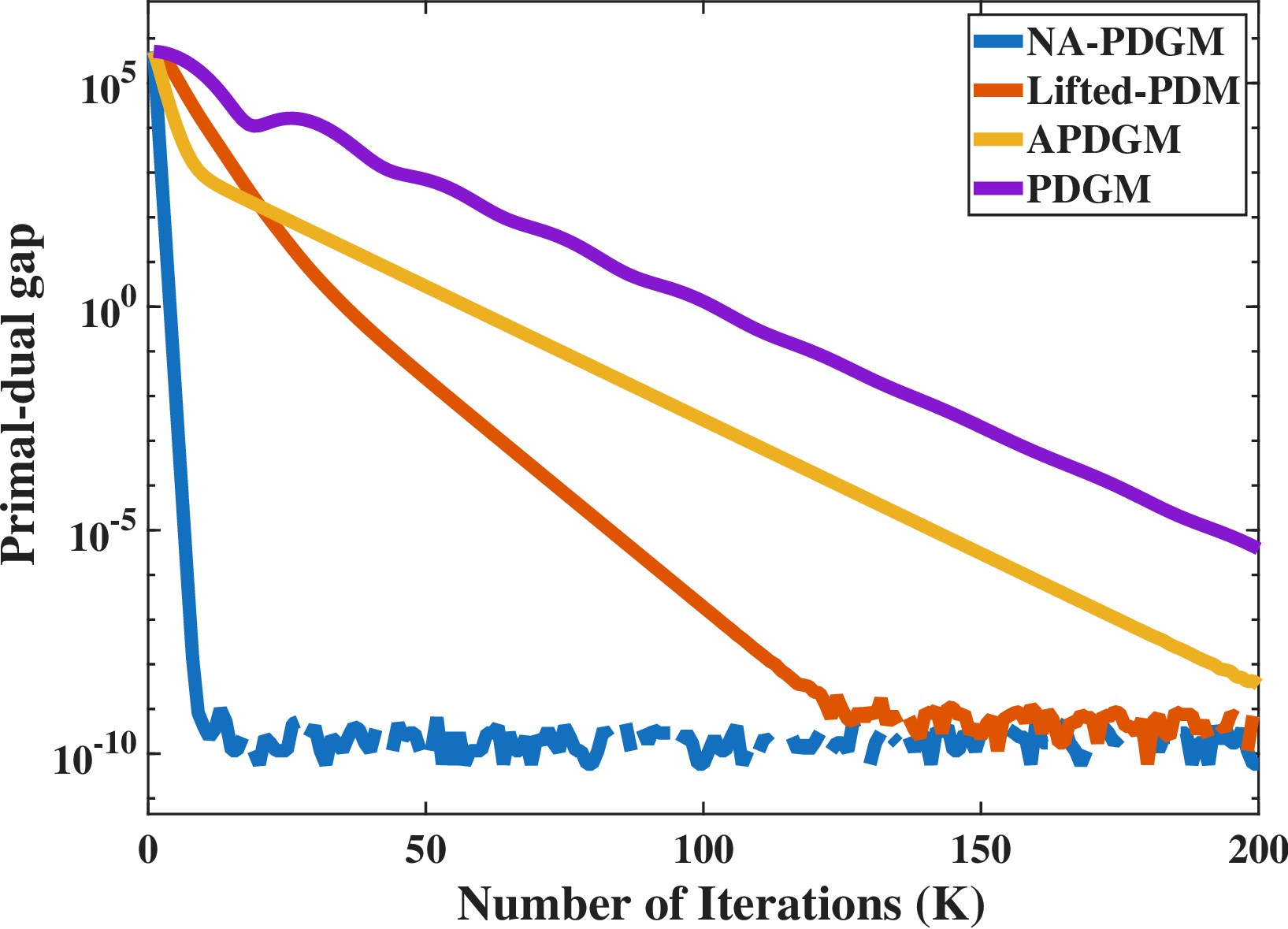}

\vspace{0.2cm}
\centerline{Case 3: $u=1$ and $r=0.1$ ($\theta=0.8944$)}
\end{minipage}
\caption{Numerical results of algorithms for Example 2}
\label{fig2}
\end{figure}

\begin{figure}[htp]
\centering
\begin{minipage}{\linewidth}
\centering
\includegraphics[width=0.48\linewidth]{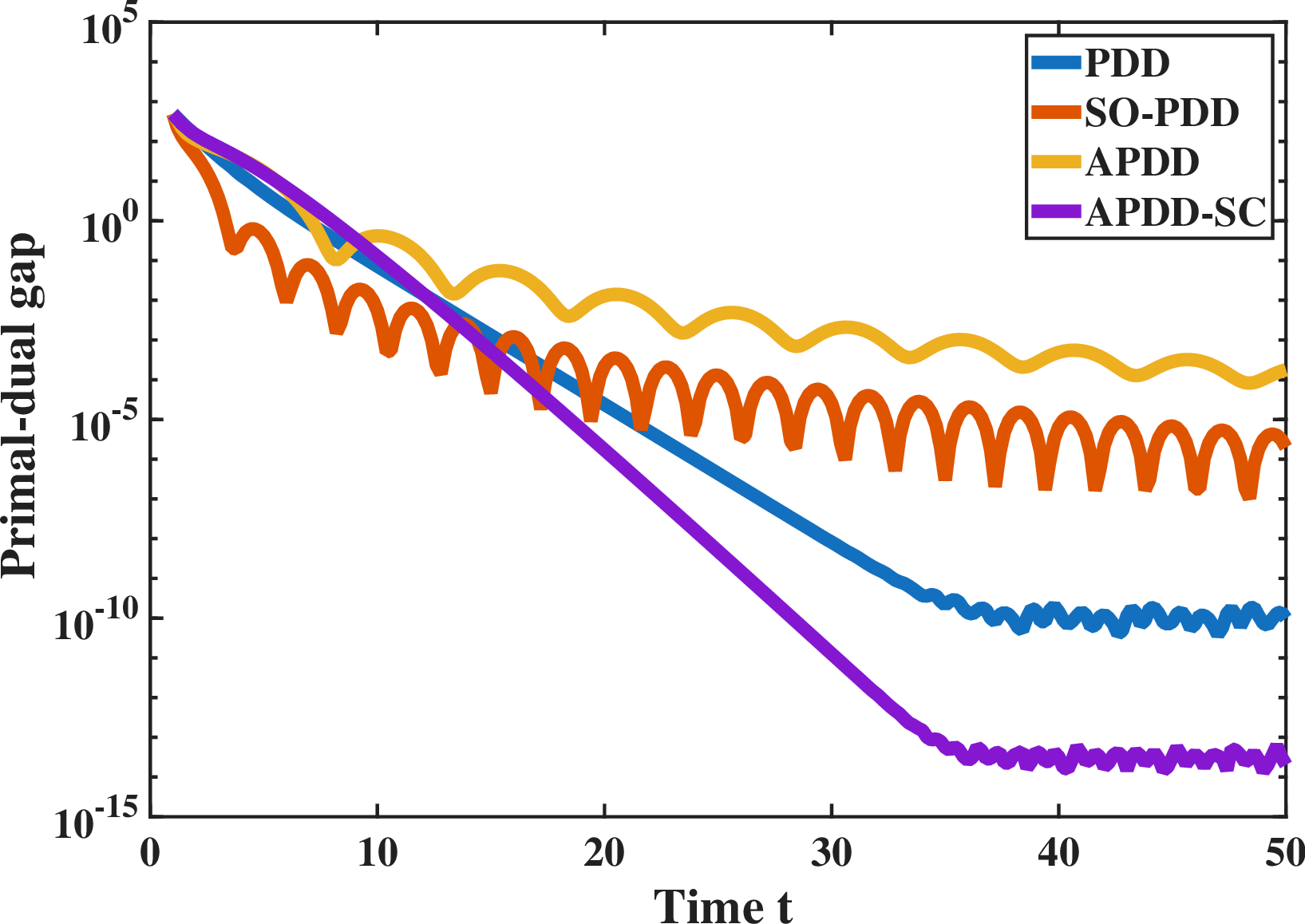}\hfill
\includegraphics[width=0.48\linewidth]{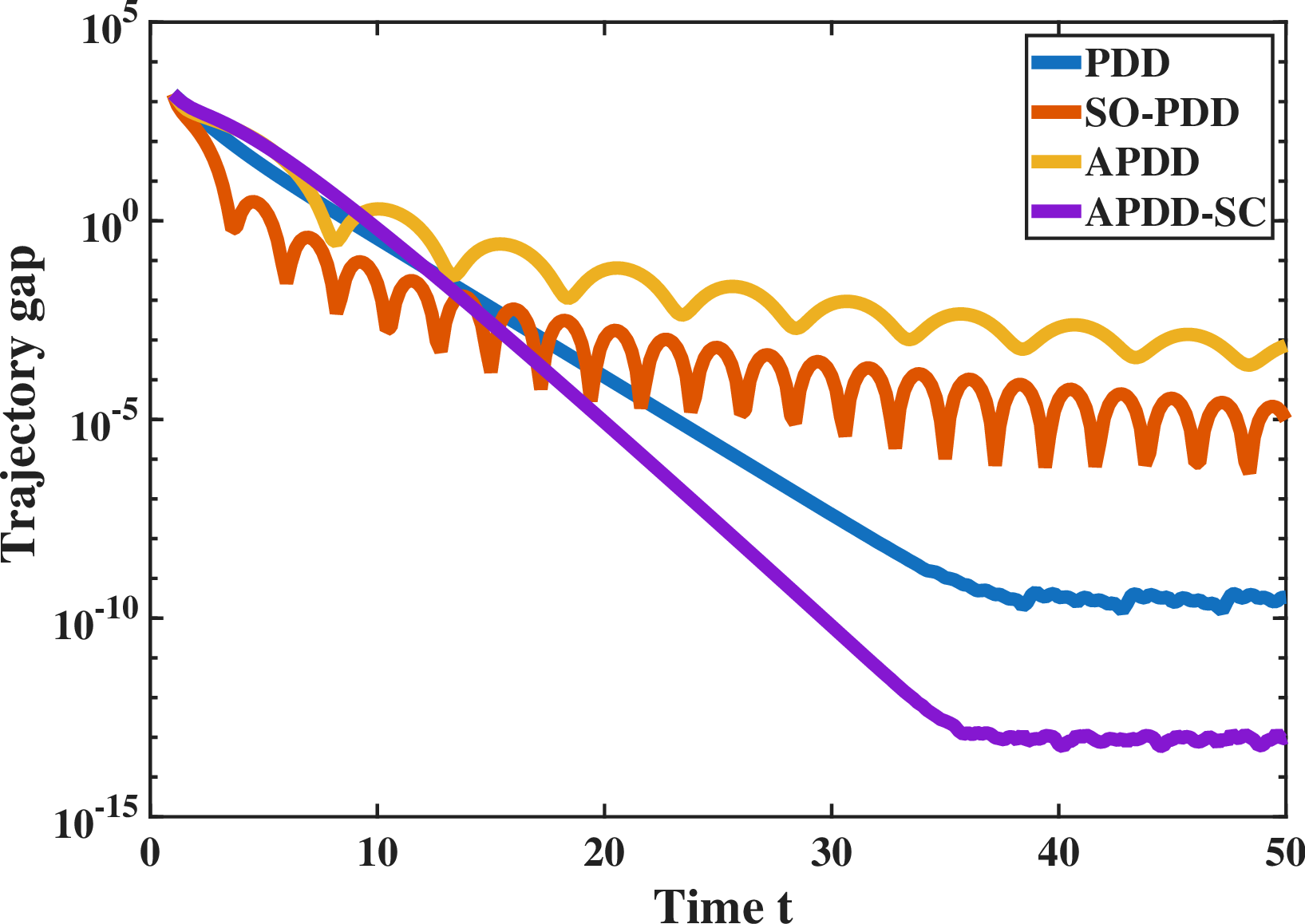}
\end{minipage}
\caption{Convergence properties  of dynamics for Example 3}  \label{fig3}
\centering
\end{figure}

\emph{Example 2.} Let $\{(a_i, b_i)\}_{i=1}^N$ denote the training dataset, where $a_i \in \mathbb{R}^n$ represents the feature vector and $b_i \in \{-1,1\}$ is the class label. We consider a regularized constrained classification problem formulated as problem~\eqref{ques}, where $f$ denotes the primal loss, $g$ the dual regularization term, and $A=[a_1,\ldots,a_N]^T\in\mathbb{R}^{N\times n}$. In this example, the primal function $f(x)$ is chosen as the logistic loss with an additional $\ell_2$ regularization term:
\[
f(x) = \tfrac{1}{N} \sum_{i=1}^{N}
\log\big(1 + \exp(-b_i a_i^T x)\big) + \tfrac{\mu}{2}\|x\|^2,
\]
which is $\mu$-strongly convex and $L$-smooth with $L=\max\{{\|a_i\|^2}/{4}+\mu|i=1,2,\cdots, N\}$.
The first term in $f$ measures the average classification error over the samples with labels $b_i \in \{-1,1\}$, while the second term penalizes the model complexity. For the dual function, we adopt a quadratic regularization $g(y) = \tfrac{r}{2}\|y\|^2+\langle y,c\rangle$ with $c\in\mathbb{R}^{N}$, which is both $r$-strongly convex and $r$-smooth. This formulation captures the trade-off between minimizing the prediction loss and enforcing the linear constraint $Ax=c$ through dual penalization. In numerical experiments, we set $N=200$ and $p=5000$, and generate $(A,b,c)$ randomly. 
We compare the performance of the proposed NA-PDGM algorithm (set $r=1/L$, $s=1/r$) with the Lifted-PDM algorithm~\citep{Thekumparampil2022}, the APDGM algorithm~\citep{Kovalev2022}, and the classical primal-dual gradient method (PDGM)~\citep{DuH19}, which does not incorporate momentum or extrapolation terms. 
As shown in Fig.~\ref{fig2}, under three different problem parameter settings and for different choices of $\theta$ in Algorithm 1, the proposed NA-PDGM algorithm consistently exhibits superior convergence behavior, which further illustrates the effectiveness of the extrapolation and momentum mechanisms.

\emph{Example 3.} Consider the $\ell_2$-regularized problem:
\begin{equation*}
	\min_{x\in\mathbb{R}^n} \frac{1}{2}\|Kx-b\|^2+\frac{\mu}{2}\|x\|^2,
\end{equation*}
where $K\in\mathbb{R}^{m\times n}$ and $b\in \mathbb{R}^{m}$. This can be rewritten as the following saddle point problem:
\begin{equation*}
	\min_{x\in\mathbb{R}^n}\max_{y\in\mathbb{R}^m} \frac{\mu}{2}\|x\|^2+\langle Kx,y\rangle-(\frac{1}{2}\|y\|^2+\langle b,y\rangle).
\end{equation*}
 In the numerical experiments, we set $m = 500$, $n= 1000$ and $\mu=0.4$, with $K$ and $ b$ randomly generated from a standard Gaussian distribution. We solve problem \eqref{ques_quad} using our accelerated primal-dual dynamic \eqref{dynamic} (APDD-SC), the unaccelerated first-order primal-dual dynamic \eqref{dy_fd} (PDD) in \citep{CherukuriSICON}, the second-order primal-dual dynamic (SO-PDD)  in \citet{HeAMO} with $\alpha(t)=5/t$, $\beta(t)=5 $, $\delta = t/3$, and the accelerated primal-dual dynamic (APDD) in \citet{ZengIFAC} with $\alpha=5$. Fig. \ref{fig3} illustrates the convergence properties of the primal-dual gap $\mathcal{L}(x(t),y^*)-\mathcal{L}(x^*,y(t))$ and the trajectory gap $\|x(t)-x^*\|^2+\|y(t)-y^*\|^2$ under different dynamics. The results show that our dynamical system \eqref{dynamic} is both faster and more stable than the second-order primal-dual dynamic in \citet{HeAMO} and the accelerated primal-dual dynamic in \citet{ZengIFAC}, and shows acceleration compare with first-order primal-dual dynamic \eqref{dy_fd}. The numerical results further demonstrate that our dynamical system enjoys linear convergence,  which is perfectly aligned with the theoretical convergence rate.

\section{Conclusion and perspective}\label{sec5}

In this paper, we study saddle point problems with bilinear coupling in the strongly convex-strongly concave setting and propose an accelerated primal-dual framework that unifies discrete and continuous-time perspectives. Specifically, we develop a first-order primal-dual gradient method that incorporates Nesterov acceleration, achieving the linear convergence rate \eqref{opt_conv} for both the primal-dual gap and the iterative gap. We further extend this idea to a continuous-time primal-dual dynamical system with constant damping, for which we establish the linear convergence rate of $\mathcal{O}(e^{-\min{\sqrt{\mu_F}, \sqrt{\mu_G}}t})$. Notably, when the bilinear term $\langle Ax, y\rangle$ vanishes, both the discrete and continuous methods naturally reduce to Nesterov accelerated schemes for unconstrained strongly convex optimization, demonstrating the generality and theoretical consistency of our approach. Numerical experiments confirm the effectiveness of the proposed methods, validating the theoretical guarantees.

 The present work is restricted to the centralized and full-information setting, which provides a basic model for studying acceleration in coupled primal-dual dynamics. Several directions deserve further investigation. One is to extend the proposed framework to distributed or partial-information scenarios, where the interaction among agents and the use of estimation or communication mechanisms may substantially affect the design and analysis of accelerated schemes. Another direction is to consider more general saddle-point models with non-bilinear coupling structures, for which the current analysis may no longer apply directly. It is also of interest to examine whether the present approach can be extended to broader problem classes, including nonconvex-nonconcave formulations arising in game-theoretic and learning-based applications.

%
%
%

\end{document}